\documentclass[a4paper,11pt]{article}
\usepackage{amsmath}
\usepackage{amsfonts}
\usepackage{amssymb}
\newtheorem{theorem}{Theorem}[section]
\newtheorem{lemma}{Lemma}[section]

\newtheorem{definition}{Definition}
\newtheorem{proposition}{Proposition}[section]
\newenvironment{proof}{{\bf Proof.}}{\par\hspace{25em}\rule{1ex}{1ex}\par}
\newcommand{\la}{\big<}
\newcommand{\ra}{\big>}
\newcommand{\ds}{\displaystyle}
\newcommand{\vf}{\vphantom{\Big|}}
\newcommand{\bml}{\begin{multline*}}
\newcommand{\eml}{\end{multline*}}
\newcommand{\ve}{\varepsilon}
\newcommand{\te}{\,\tilde e}
\newcommand{\tf}{\,\tilde f}
\newcommand{\bt}{\beta\,}
\newcommand{\al}{\alpha\,}
\newcommand{\dt}{\delta\,}

\title{ Invariant totally geodesic unit vector fields on three-dimensional Lie groups.}
\author{Yampolsky A.}
\date{}
\begin{document}
\maketitle
\begin{abstract}
We give a complete list of those left invariant unit vector fields on three-dimensional Lie groups
with the left-invariant metric that generate a totally geodesic submanifold in the unit tangent
bundle of a group with the Sasaki metric. As a result, each class of three-dimensional Lie groups
admits the totally geodesic unit vector field. From geometrical viewpoint, the field is either
parallel or characteristic vector field of a natural almost contact structure on the group.


{\it AMS subject class:} Primary 53B20, 53B25; Secondary 53C25.

\end{abstract}
\section*{Introduction}
Let $(M^n,g)$ be Riemannian manifold and $(T_1M^n, g_s)$ its unit tangent bundle with Sasaki
metric. Consider a unit vector field $\xi$ as a (local) mapping
$$
\xi:M^n\to T_1M^n.
$$
\begin{definition}\label{def:1}
A unit vector field  $\xi$  on Riemannian manifold $M^n$ is called totally geodesic if the image of
(local) imbedding $\xi:M^n\to T_1M^n$ is totally geodesic submanifold in the unit tangent bundle
$T_1M^n$ with Sasaki metric.
\end{definition}
In a similar way one can define a \emph{locally minimal} unit vector field as the field of
\emph{zero mean curvature}. A number of examples of locally minimal unit vector fields was found
recently by L.~Vanhecke, E.~Boeckx, K.~Tsukada, J.C.~Gonz\'alez -D\'avila, O.~Gil-Medrano and
others \cite{BX-V1, BX-V2, BX-V4, B-Ch-N, GM, GM-GD-Vh, GM-LF, GD-V1, GD-V2, TS-V1,TS-V2,TS-V3}.
Particulary, K~Tsukada and L.~Vanhecke \cite{TS-V2} described all minimal left-invariant unit
vector fields on three-dimensional Lie groups with the left-invariant metric.

The key step to the totally geodesic unit vector fields was made in \cite{Ym1}, where the author
have found the second fundamental form of $\xi(M^n)$ explicitly, using a special normal frame. This
expression allowed, also, to find an examples of unit vector fields of \emph{constant mean
curvature}. Using this expression, the author described all the 2-manifolds that admit a totally
geodesic unit vector field and the field itself \cite{Ym5}. In the case of higher dimensions only
partial results are known. The most general states that if $M^{2m+1}$ is a Sasakian manifold and
$\xi$ is a characteristic vector field of the Sasakian structure, then $\xi(M^{2m+1})$ is totally
geodesic in $T_1M^{2m+1}$ \cite{Ym3}.

Particularly, the Hopf unit vector field on a unit $S^{2m+1}$ is totally geodesic. More
specifically, the Hopf vector field belongs to the class of left invariant unit vector fields on
$S^3$ as a Lie group with the left-invariant Riemannian metric. In this paper, we give full
description of 3-dimensional Lie groups with left-invariant metric which admit a totally geodesic
left-invariant unit vector fields and the fields themselves. As a consequence we have found that,
in non-trivial case, for each totally geodesic left invariant  unit vector field $\xi$ the
structure $\big(\phi=-\nabla\xi,\ \xi, \ \eta=\la\xi,\cdot\ra\big)$ is an almost contact one on the
corresponding Lie group and $\xi$  is a characteristic vector field of this structure. If $\xi$ is
a Killing unit vector field, then the structure is Sasakian.

The paper is organized as follows. In Section 1 we give some preliminaries. In Section 2 we
consider the unimodular Lie groups. We prove that if the totally geodesic unit vector field exists
on a given group, then it is an eigenvector of the Ricci tensor which corresponds to the Ricci
principal curvature $\rho=2$ (Theorem \ref{Thrm:2.1}). The Theorem \ref{Thrm:2.2} provides the
complete list of totally geodesic unit vector fields on a corresponding Lie group as well as the
conditions on structure constants of the group. In a series of Propositions 2.2 -- 2.6, we give a
description of totally geodesic unit vector field in unimodular case from the contact geometry
viewpoint.

In Section 3 we consider the non-unimodular case. The Theorem  \ref{Thrm:3.2} provides an explicit
expression for the totally geodesic unit vector field as well as the conditions on structure
constants of the corresponding group. Finally, the Proposition \ref{prop:3.1} gives the geometrical
characterization of the totally geodesic unit vector field and clarifies the structure of the
corresponding non-unimodular Lie group.

\section{Some preliminaries}
Let $(M^n,g)$ be Riemannian manifold. Denote by $\nabla$ the Levi-Civita connection on $M^n$.
Introduce a pointwise linear operator $A_\xi:T_qM^n\to \xi^\perp_q$, acting as
$$
A_\xi X=-\nabla_X\xi.
$$
In case of integrable distribution $\xi^\perp$, the unit vector field $\xi$ is called
\emph{holonomic}. In this case the operator $A_\xi$ is symmetric and is known as Weingarten or a
\emph{shape operator} for each hypersurface of the foliation. In general, $A_\xi$ is not symmetric
but formally preserves the Codazzi equation. Namely, a covariant derivative of $A_\xi$ is defined
by
\begin{equation}\label{eqno:1}
(\nabla_X A_\xi) Y=-\nabla_X\nabla_Y\xi+\nabla_{\nabla_XY}\xi.
\end{equation}
Then for the curvature operator of $M^n$ we can write down the Codazzi-type equation
$$
R(X,Y)\xi=(\nabla_Y A_\xi) X-(\nabla_X A_\xi) Y.
$$
From this viewpoint, it is natural to call the operator $A_\xi$ by \emph{non-holonomic shape}
operator.

Introduce a symmetric tensor field
\begin{equation}\label{eqno:2}
Hess_\xi(X,Y)=\frac12\big[(\nabla_Y A_\xi) X+(\nabla_X A_\xi) Y\big],
\end{equation}
which is a \textit{symmetric part} of covariant derivative of $A_\xi$. The trace
$$
-\sum_{i=1}^n Hess_\xi(e_i,e_i):=\Delta \xi,
$$
where $e_1,\dots e_n$ is an orthonormal frame, is known as \textit{rough Laplacian} \cite{Besse} of
the field $\xi$. Therefore, one can treat the tensor field \eqref{eqno:2} as a \textit{rough
Hessian} of the field.

For the mapping $f:(M,g)\to (N,h)$ between Riemannian manifolds the \emph{energy } of $f$ is
defined as
$$
E(f):=\frac12 \int_M |d\,f|^2\,d\,Vol_M,
$$
where $|d\,f|$ is a norm of 1-form $d\,f$ in the cotangent bundle $T^*M$. The mapping $f$ is called
harmonic if it is a critical point of the functional $E(f)$. Supposing on $T_1M$ the Sasaki metric,
a unit vector field is called \textit{harmonic}, if it is a critical point of energy functional of
mapping $\xi:M^n\to T_1M^n$. This definition presumes the variation within the class of unit vector
fields.  From this viewpoint, the unit vector field is harmonic if and only if \cite{Weigm}
$$
\Delta\xi=-|\nabla\xi|^2\xi.
$$
There exist the unit vector fields that fail to be critical within a wider class of all mappings
$f:M^n\to T_1M^n$ \cite{GM}.  Introduce a tensor field
\begin{equation}\label{eqno:3}
Hm_\xi(X,Y)=\frac12\big[R(\xi,A_\xi X)Y+R(\xi,A_\xi Y)X\big].
\end{equation}
A \emph{harmonic} unit vector field $\xi$ defines a \emph{harmonic mapping} $\xi:M^n\to T_1M^n$ if
and only if  \cite{GM}
$$
\sum_{i=1}^n\, Hm_\xi(e_i,e_i)=0.
$$
From this viewpoint, it is natural to call the tensor field \eqref{eqno:3} by \textit{harmonicity
tensor} of the field $\xi$.

In terms of the tensors $Hess_\xi$ and $Hm_\xi$ the conditions on $\xi$ to be totally geodesic are
as follows \cite{Ym6}.
\begin{theorem}\label{theo:1} A unit vector  field $\xi$ on a given Riemannian manifold $M^n$ is
totally geodesic if and only if
$$
Hess_\xi(X,Y)+A_\xi Hm_\xi(X,Y)-\la A_\xi X, A_\xi Y\ra \xi=0
$$
for all vector fields $X,Y$ on $M^n$.
\end{theorem}
It is natural to introduce a tensor field
\begin{equation}\label{eqno:4}
TG_\xi(X,Y)=Hess_\xi(X,Y)+A_\xi Hm_\xi(X,Y)-\la A_\xi X, A_\xi Y\ra \xi
\end{equation}
as a \emph{total geodesity} tensor field.

The treatment of 3-dimensional Lie groups is based on J.~Milnor description of 3-dimensional Lie
groups via the structure constants \cite{Mn}.

In the case of  \emph{unimodular} Lie group with the left-invariant metric, there is an orthonormal
frame $e_1,e_2,e_3$ of its Lie algebra such that the bracket operations are defined by
\begin{equation}\label{eqno:5}
[e_2,e_3]=\lambda_1e_1,\quad [e_3,e_1]=\lambda_2e_2,\quad[e_1,e_2]=\lambda_3e_3.
\end{equation}
The constants $\lambda_1, \lambda_2, \lambda_3$ completely determine the topological structure of
corresponding Lie group as in the following table:
$$
\begin{tabular}{|c|c|}
  \hline
  signs of $\lambda_1, \lambda_2, \lambda_3$ & Associated Lie grous \\
  \hline
  $+,+,+$ & $SU(2)$ or $SO(3)$ \\
  $+,+,-$ & $SL(2,\mathbb{R})$ or $O(1,2)$ \\
  $+,+,0$ & $E(2)$ \\
  $+,-,0$ & $E(1,1)$\\
  $+,0,0$ & $Nil^3$ (Heisenberg group)  \\
  $0,0,0$ & $\mathbb{R}\oplus\mathbb{R}\oplus \mathbb{R}$\\
  \hline
\end{tabular}
$$
In the case of non-unimodular Lie group, let $e_1$ be a unit vector orthogonal to the unimodular
kernel $U$ and choose an orthonormal basis $\{e_2,e_3\}$ of $U$ which diagonalizes the symmetric
part of $ad_{e_1}\big|_U$. Then the bracket operation can be expressed as
\begin{equation}\label{eqno:6}
[e_1,e_2]=\al e_2+\bt e_3, \quad [e_1,e_3]=-\bt e_2+\dt   e_3, \quad  [e_2,e_3]=0.
\end{equation}
If necessary, changing $e_1$ to $-e_1$, we can assume $\al+\dt  >0$ and by possibly alternating
$e_2$ and $e_3$, we may also suppose $\al\geq \dt$ \cite{TS-V2}.

\section{The unimodular case}

Choose the orthonormal frame as in \eqref{eqno:5}. Define a connection numbers by
$$
\mu_i=\frac12(\lambda_1+\lambda_2+\lambda_3)-\lambda_i.
$$
Then the Levi-Civita covariant derivatives can be expressed via the cross-products as follows
\begin{equation}\label{eqno:7}
\nabla_{e_i}e_k=\mu_i\,e_i\times e_k.
\end{equation}
For any left-invariant unit vector field $\xi=x_1e_1+x_2e_2+x_3e_3$ we have
\begin{equation}\label{eqno:8}
\nabla_{e_i} \xi=\mu_i\, e_i\times \xi.
\end{equation}
Denote $ N_i=e_i\times \xi$.  Then
\begin{equation}\label{eqno:9}
\nabla_{e_i}\xi=\mu_i\, e_i\times \xi=\mu_i\, N_i.
\end{equation}
As a consequence, the matrix of the Weingarten operator takes the form
\begin{equation}\label{W}
A_\xi= \left(
  \begin{array}{ccc}
    0 & -\mu_2 x_3 & \mu_3 x_2 \\
    \mu_1 x_3 & 0 & -\mu_3 x_1 \\
    -\mu_1 x_2 & \mu_2 x_1 & 0 \\
  \end{array}
\right)
\end{equation}
We will need the following technical Lemma.

\begin{lemma}\label{Lemma:2.1}
Let $G$ be a three-dimensional unimodular Lie group with the left-invariant metric and let
$\{e_i,\, i=1,2,3\}$ be an orthonormal basis for the Lie algebra satisfying \eqref{eqno:5}.  Then
for any left-invariant unit vector field $\xi=x_1e_1+x_2e_2+x_3e_3$ we have
$$
\begin{array}{l}
A_\xi e_i=-\mu_i\, e_i\times\xi=-\mu_i\, N_i\\[1ex]
 (\nabla_{e_i}A_\xi)e_i=\mu_i^2(\xi-x_ie_i),\\[1ex]
 (\nabla_{e_i}A_\xi)e_k=\ve_{ikm}\mu_i\mu_m N_m-\mu_i\mu_k x_ie_k \quad (i\ne k),\\[1ex]
 R(e_i,e_k)\xi=-\ve_{ikm}\sigma_{ik}N_m,
 \end{array}
$$
 where $\sigma_{ik}=\sigma_{ki}=\mu_i\mu_m+\mu_k\mu_m-\mu_i\mu_k$ and  $\ve_{ikm}=\la e_i\times e_k, e_m\ra$.
\end{lemma}
\begin{proof}
The first equality comes from definitions. For the rest, we have
$$
\begin{array}{l}
 \nabla_{e_i}e_k=\mu_i\,e_i\times e_k=\ve_{ikm}\mu_i\,e_m,\quad  \nabla_{\nabla{e_i}e_k}\xi=\ve_{ikm}\mu_i\mu_mN_m,\\[1ex]
 \nabla_{e_i}\nabla_{e_k}\xi=\mu_i\mu_k\, e_i\times (e_k \times \xi)=\mu_i\mu_k(x_ie_k-\delta_{ik}\xi).
\end{array}
$$
Therefore,
$$
(\nabla_{e_i}A_\xi)e_k= \nabla_{\nabla{e_i}e_k}\xi-
\nabla_{e_i}\nabla_{e_k}\xi=\ve_{ikm}\mu_i\mu_mN_m+\mu_i\mu_k(\delta_{ik}\xi-x_ie_k)
$$
Setting $i=k$ and then $i\ne k$, we get the second and the third equalities. From Codazzi equation
$$
\begin{array}{ll}
R(e_i,e_k)\xi=&(\nabla_{e_k}A_\xi)e_i-(\nabla_{e_i}A_\xi)e_k=\\[1ex]
&\ve_{kim}\mu_k\mu_mN_m+\mu_k\mu_i(\delta_{ki}\xi-x_ke_i)-\\
&\ve_{ikm}\mu_i\mu_mN_m-\mu_i\mu_k(\delta_{ik}\xi-x_ie_k)=\\[1ex]
&-\ve_{ikm}(\mu_i\mu_m+\mu_k\mu_m)N_m+\mu_i\mu_k(x_ie_k-x_ke_i).
\end{array}
$$
Remark, that $N_m=\ve_{ikm}(x_ie_k-x_ke_i)$ and hence

$$
 R(e_i,e_k)\xi= -\ve_{ikm}(\mu_i\mu_m+\mu_k\mu_m-\mu_i\mu_k)N_m
 $$
\end{proof}
Remark that chosen frame diagonalises the Ricci tensor \cite{Mn}. Moreover,
$$
2\mu_i\mu_k=\rho_m,
$$
where $\rho _m$ is the principal Ricci curvature and $i\ne k\ne m$. It also worthwhile to mention
that
$$
\sigma_{ik}=\frac12(\rho_k+\rho_i-\rho_m)
$$
is nothing else bur the sectional curvature of the left-invariant metric in a direction of
$e_i\wedge e_k$.
\begin{lemma}\label{Lemma:2.2}
Let $G$ be a three-dimensional unimodular Lie group with the left-invariant metric and let
$\{e_i,\, i=1,2,3\}$  be an orthonormal basis for the Lie algebra satisfying \eqref{eqno:5}. Then a
left-invariant unit vector field $\xi=x_1e_1+x_2e_2+x_3e_3$ is totally geodesic if and only if for
any $i\ne k\ne m$
$$
\begin{array}{ll}
TG(e_i,e_i)=& x_i\mu_i\Big\{x_m(\sigma_{ik}\mu_k-\mu_i)N_k-x_k(\sigma_{im}\mu_m-\mu_i)N_m\Big\}=0,\\[2ex]
2TG(e_i,e_k)=&\ve_{ikm}\Big\{ -x_ix_m\mu_i(\sigma_{ik}\mu_i-\mu_k)N_i+x_kx_m\mu_k(\sigma_{ik}\mu_k-\mu_i)N_k+\\[1ex]
&\Big(\mu_i\mu_m(1-\sigma_{km})-\mu_k\mu_m(1-\sigma_{im})+\mu_i(\sigma_{km}\mu_m-\mu_k)x_i^2-\\
&\hspace{4.7cm}\mu_k(\sigma_{im}\mu_m-\mu_i)x_k^2\Big)N_m\Big\}=0,
\end{array}
$$
where $\sigma_{ik}=\sigma_{ki}=\mu_i\mu_m+\mu_k\mu_m-\mu_i\mu_k$ and  $\ve_{ikm}=\la e_i\times e_k,
e_m\ra$.
\end{lemma}
\begin{proof}
Calculate $Hess_\xi(e_i,e_i)-|A_\xi e_i|^2\,\xi$. We have
\begin{equation}
\begin{array}{ll}\label{eqno:11}
(\nabla_{e_i}A_\xi)e_i-|A_\xi e_i|^2\xi=&\mu_i^2(\xi-x_ie_i)-\mu_i^2(1-x_i^2)\xi=-\mu_i^2x_i(e_i-x_i\xi)=\\
&-\mu_i^2x_i((1-x_i^2)e_i-x_ix_ke_k-x_ix_me_m)=\\
&-\mu_i^2x_i((x_k^2+x_m^2)e_i-x_ix_ke_k-x_ix_me_m)=\\
&-\mu_i^2x_i(x_k(x_ke_i-x_ie_k)+x_m(x_me_i-x_ie_m)=\\
&-\mu_i^2x_i(x_k\ve_{kim}N_m+x_m\ve_{mik}N_k)=\\[1ex]
&\ve_{ikm}\mu_i^2x_i(x_kN_m-x_mN_k).
\end{array}
\end{equation}
Find now $A_\xi Hm_\xi(e_i,e_i)$.  Using Lemma \ref{Lemma:2.1}, we have
$$
\begin{array}{l}
Hm_\xi(e_i,e_i)=R(\xi,A_\xi e_i)e_i=\\[1ex]
\hspace{2.5cm}\la R(\xi,A_\xi e_i)e_i,e_k\ra e_k+\la R(\xi,A_\xi e_i)e_i,e_m\ra e_m=\\[1ex]
\hspace{2.5cm} \la R(e_i,e_k)\xi,A_\xi e_i\ra e_k+\la R(e_i,e_m)\xi,A_\xi e_i \ra e_m=\\
\hspace{1.5cm}\mu_i\Big(\ve_{ikm}\sigma_{ik}\la e_m\times\xi,e_i\times \xi\ra e_k+
\ve_{imk}\sigma_{im}\la e_k\times\xi,e_i\times \xi,\ra e_m \Big)=\\
\hspace{2.5cm}-\mu_i\ve_{ikm}\Big(x_ix_m\sigma_{ik}e_k-x_ix_k\sigma_{im}e_m\Big)
\end{array}
$$
Therefore,
\begin{equation}\label{eqno:12}
A_\xi Hm_\xi(e_i,e_i)=\ve_{ikm}\mu_ix_i\Big(x_m\sigma_{ik}\mu_kN_k-x_k\sigma_{im}\mu_mN_m\Big)
\end{equation}
Adding \eqref{eqno:11} and \eqref{eqno:12}, after evident simplifications we get $TG_\xi(e_i,e_i)$.

Applying Lemma \ref{Lemma:2.1} for $i\ne k$, we get
$$
\begin{array}{ll}
2Hess_\xi(e_i,e_k)=&(\nabla_{e_i}A_\xi)e_k+(\nabla_{e_k}A_\xi)e_i=\\[1ex]
&\ve_{ikm}(\mu_i\mu_m-\mu_k\mu_m)N_m-\mu_i\mu_k(x_ie_k+x_ke_i).
\end{array}
$$
Evidently,
$$
\la A_\xi e_i,A_\xi e_k\ra\xi=\mu_i\mu_k\la e_i\times \xi,e_k\times\xi\ra\xi=-\mu_i\mu_k x_i x_k
\xi
$$
Subtracting, we get
\begin{multline*}
2Hess_\xi(e_i,e_k)-2\la A_\xi e_i,A_\xi
e_k\ra\xi=\ve_{ikm}(\mu_i\mu_m-\mu_k\mu_m)N_m-\\
\mu_i\mu_k(x_ie_k+x_ke_i-2x_ix_k\xi).
\end{multline*}
Observe that
$$
\begin{array}{l}
x_ie_k+x_ke_i-2x_ix_k\xi=x_i(1-2x_k^2)e_k+x_k(1-2x_i^2)e_i-2x_ix_kx_me_m=\\[1ex]
\hspace{2cm} x_k(x_k^2-x_i^2+x_m^2)e_i+ x_i(-x_k^2+x_i^2+x_m^2)e_k-2x_ix_kx_me_m=\\[1ex]
\hspace{0.5cm} x_kx_m(x_me_i-x_ie_m)+x_ix_m(x_me_k-x_ke_m)+(x_k^2-x_i^2)(x_ke_i-x_ie_k)=\\[1ex]
\hspace{2cm}x_kx_m\ve_{mik}e_k\times\xi+x_ix_m\ve_{mki}e_i\times\xi+(x_k^2-x_i^2)\ve_{kim}e_m\times\xi=\\[1ex]
\hspace{2cm}\ve_{ikm}\big(x_kx_mN_k-x_ix_m N_i-(x_k^2-x_i^2)N_m\big).
\end{array}
$$
Therefore,
\begin{multline*}
2Hess_\xi(e_i,e_k)-2\la A_\xi e_i,A_\xi e_k\ra\xi=\\
\ve_{ikm}\Big\{\mu_i\mu_kx_m(-x_kN_k+x_iN_i)+(\mu_i\mu_m-\mu_k\mu_m+(x_k^2-x_i^2)\mu_i\mu_k)N_m\Big\}
\end{multline*}
To find $Hm_\xi(e_i,e_k)$, calculate  $R(\xi,A_\xi e_i)e_k$.  We have
$$
\begin{array}{ll}
R(\xi,A_\xi e_i)e_k=&\la R(\xi,A_\xi e_i)e_k,e_i\ra e_i+\la R(\xi,A_\xi e_i)e_k,e_m\ra e_m=\\
&\la R(e_k,e_i)\xi,A_\xi e_i\ra e_i+ \la R(e_k,e_m)\xi,A_\xi e_i\ra e_m=\\
&\mu_i\sigma_{ki}\ve_{kim}\la e_m\times\xi,e_i\times
\xi\ra e_i+\mu_i\sigma_{km}\ve_{kmi}\la e_i\times\xi,e_i\times \xi\ra e_m=\\[1ex]
&\ve_{ikm}\big\{\mu_i\sigma_{ki}x_ix_me_i+\mu_i\sigma_{km}(1-x_i^2)e_m\big\}
\end{array}
$$
Therefore,
$$
\begin{array}{ll}
A_\xi R(\xi,A_\xi e_i)e_k=&
-\ve_{ikm}\big\{\mu_i^2\sigma_{ki}x_ix_me_i\times\xi+\mu_i\mu_m\sigma_{km}(1-x_i^2)e_m\times\xi\big\}=\\[1ex]
&-\ve_{ikm}\big\{\mu_i^2\sigma_{ki}x_ix_mN_i+\mu_i\mu_m\sigma_{km}(1-x_i^2)N_m\big\}
\end{array}
$$
Thus,
\begin{multline*}
2A_\xi Hm_\xi(e_i,e_k)= \ve_{ikm}\big\{-\mu_i^2\sigma_{ki}x_i x_m
N_i+\mu_k^2\sigma_{ki}x_kx_mN_k-\\
\big(\mu_i\mu_m\sigma_{km}(1-x_i^2)-\mu_k\mu_m\sigma_{im}(1-x_k^2)\big)N_m\big\}
\end{multline*}
So, finally
$$
\begin{array}{l}
2\ve_{ikm}TG_\xi(e_i,e_k)=
x_ix_m(-\sigma_{ik} \mu_i^2+\mu_i\mu_k)N_i-x_kx_m(-\sigma_{ik} \mu_k^2+\mu_i\mu_k)N_k+ \\
\Big(\mu_i\mu_m-\mu_k\mu_m-\mu_i\mu_m\sigma_{km}(1-x_i^2)+\mu_k\mu_m\sigma_{im}(1-x_k^2)+\\
\hspace{8cm} \mu_i\mu_k(x_k^2-x_i^2)\Big)N_m=\\
x_ix_m\mu_i(-\sigma_{ik} \mu_i+\mu_k)N_i-x_kx_m\mu_k(-\sigma_{ik} \mu_k+\mu_i)N_k+ \\
\Big(\mu_i\mu_m(1-\sigma_{km})-\mu_k\mu_m(1-\sigma_{im})+\mu_i(\sigma_{km}\mu_m-\mu_k)x_i^2-\\
\hspace{8cm}\mu_k(\sigma_{im}\mu_m-\mu_i)x_k^2\Big)N_m.
\end{array}
$$
The proof is complete.
\end{proof}

\begin{theorem}\label{Thrm:2.1} Let $G$ be a three-dimensional unimodular Lie group with the left-invariant metric and let $\{e_i,\,
i=1,2,3\}$  be an orthonormal basis for the Lie algebra satisfying \eqref{eqno:5}. Denote by
$\rho_1,\rho_2,\rho_3$ the principal Ricci curvatures of the given group. Then the set of
left-invariant totally geodesic unit vector fields can be described as follows.
\begin{center}
Table 1
\end{center}
\begin{center}
\begin{tabular}{|c|c|c|r|r|r||c|}
  \hline
 \vphantom{\Big|}$\rho_1$ & $\rho_2$ & $\rho_3$ & $\mu_1$ & $\mu_2$ & $\mu_3$ & $\xi$ \\
  \hline\hline
 \vphantom{\Big|} 0 & 0 & 0 & 0 & 0 & 0 & S \\
 \hline
\vphantom{\Big|}  0 & 0 & 0 & $\ne0$ & 0 & 0 & $\pm e_1$, $S\cap\{e_2,e_3\}_R$ \\
\hline
\vphantom{\Big|}  0 & 0 & 0 & $0$ & $\ne$ 0 & 0 & $\pm e_2$, $S\cap\{e_1,e_3\}_R$ \\
\hline
\vphantom{\Big|}  0 & 0 & 0 & $0$ & 0 & $\ne0$ & $\pm e_3$, $S\cap\{e_1,e_2\}_R$ \\
\hline \vphantom{\Big|}  2 &   &   &  &  &  & $\pm e_1$ \\
\hline
\vphantom{\Big|}    & 2 &   &  &  &  & $\pm e_2$ \\
\hline
\vphantom{\Big|}    &   & 2 &   &   &   & $\pm e_3$ \\
\hline
\vphantom{\Big|}  2 & 2 &   &   &   &   & $S\cap\{e_1,e_2\}_R$ \\
\hline
\vphantom{\Big|}  2 &   & 2 &   &   &  & $S\cap\{e_1,e_3\}_R$  \\
\hline
\vphantom{\Big|}    & 2 & 2 &   &   &   & $S\cap\{e_2,e_3\}_R$  \\
\hline
\vphantom{\Big|}  2 & 2 & 2 &   &   &  & S  \\
  \hline
\end{tabular}
\end{center}
where $S\cap\{e_i,e_k\}_R$ means the set of unit vectors in a plane, spanned by $e_i$ and $e_k$, in
the tangent space of the group at the unit element.
\end{theorem}
\begin{proof}
Rewrite the result of Lemma \ref{Lemma:2.2} for various combinations of indices.
$$
(1,1)\quad x_1\mu_1\Big\{x_3(\sigma_{12}\mu_2-\mu_1)N_2-x_2(\sigma_{13}\mu_3-\mu_1)N_3\Big\}=0,
$$
$$
(2,2)\quad x_2\mu_2\Big\{x_3(\sigma_{21}\mu_1-\mu_2)N_1-x_1(\sigma_{23}\mu_3-\mu_2)N_3\Big\}=0,
$$
$$
(3,3)\quad x_3\mu_3\Big\{x_2(\sigma_{31}\mu_1-\mu_3)N_1-x_1(\sigma_{32}\mu_2-\mu_3)N_2\Big\}=0,
$$
\begin{multline*}
(1,2)\quad -x_1x_3\mu_1(\sigma_{12}\mu_1-\mu_2)N_1+x_2x_3\mu_2(\sigma_{12}\mu_2-\mu_1)N_2+
\Big(\mu_1\mu_3(1-\sigma_{23})-\\
\mu_2\mu_3(1-\sigma_{13})+\mu_1(\sigma_{23}\mu_3-\mu_2)x_1^2-\mu_2(\sigma_{13}\mu_3-\mu_1)x_2^2\Big)N_3=0,
\end{multline*}
\begin{multline*}
(2,3)\quad -x_2x_1\mu_2(\sigma_{23}\mu_2-\mu_3)N_2+x_3x_1\mu_3(\sigma_{23}\mu_3-\mu_2)N_3+
\Big(\mu_2\mu_1(1-\sigma_{31})-\\
\mu_3\mu_1(1-\sigma_{21})+\mu_2(\sigma_{31}\mu_1-\mu_3)x_2^2-\mu_3(\sigma_{21}\mu_1-\mu_2)x_3^2\Big)N_1=0,
\end{multline*}
\begin{multline*}
(3,1)\quad -x_3x_2\mu_3(\sigma_{13}\mu_3-\mu_1)N_3+x_1x_2\mu_1(\sigma_{13}\mu_1-\mu_3)N_1+
\Big(\mu_3\mu_2(1-\sigma_{12})-\\
\mu_1\mu_2(1-\sigma_{32})+\mu_3(\sigma_{12}\mu_2-\mu_1)x_3^2-\mu_1(\sigma_{32}\mu_2-\mu_3)x_1^2\Big)N_2=0.
\end{multline*}
The vectors $N_1, N_2$ and $N_3$ are linearly dependent:
$$
x_1N_1+x_2N_2+x_3N_3=0,
$$
but linearly independent in pairs for general (not specific) field $\xi$.

\noindent\textbf{The case} $x_1\ne0,x_2\ne0,x_3\ne0$.

The subcase 1: $\mu_1=0,\mu_2=0,\mu_3=0$. All equations are fulfilled evidently. Therefore,
\emph{any left-invariant vector field is totally geodesic in this case}, and we get the first row
in the Table 1.

The subcase 2: $\mu_1=0$, $\mu_2\ne0$ or $\mu_3\ne0$. Then from (2,2) and (3,3) we see, that
$\mu_2=0,\mu_3=0$. Contradiction. In a similar way we exclude the cases when $\mu_i=0$, but
$\mu_k^2+\mu_m^2\ne 0$ for arbitrary triple of different indices $(i, k, m)$.

The subcase 3: $\mu_1\ne0,\mu_2\ne0,\mu_3\ne0$.  Since $N_1,N_2$ and $N_3$ are linearly independent
in pairs, from (1,1), (2,2) and(3,3) we conclude:
\begin{equation}\label{eqno:13}
\left\{
\begin{array}{l}
\sigma_{12}\mu_2-\mu_1=0,\\
\sigma_{12}\mu_1-\mu_2=0,
\end{array}\right.,\quad
\left\{\begin{array}{l}
\sigma_{13}\mu_3-\mu_1=0,\\
\sigma_{13}\mu_1-\mu_3=0,
\end{array}\right.
\left\{ \begin{array}{l}
\sigma_{23}\mu_2-\mu_3=0,\\
\sigma_{23}\mu_3-\mu_2=0.
\end{array}\right.
\end{equation} As a consequence,
$$
\left\{ \begin{array}{l}
(\sigma_{12}-1)(\mu_1+\mu_2)=0,\\
(\sigma_{13}-1)(\mu_1+\mu_3)=0,\\
(\sigma_{23}-1)(\mu_2+\mu_3)=0.\\
\end{array}\right.
$$
Taking  into account \eqref{eqno:13}, the rest of the equations yield
$$\left\{
\begin{array}{l}
\mu_1\mu_3(1-\sigma_{23})- \mu_2\mu_3(1-\sigma_{13})=0,\\
\mu_1\mu_2(1-\sigma_{13})- \mu_1\mu_3(1-\sigma_{12})=0,\\
\mu_2\mu_3(1-\sigma_{12})- \mu_1\mu_2(1-\sigma_{23})=0.
\end{array}\right.
$$
Since $\mu_i\ne0\ (i=1,2,3)$, we conclude $\sigma_{ik}=1 \ (i,k=1,2,3)$ and therefore $\rho_i=2 \
(i=1,2,3)$. This is the case of the last row in the Table 1.

\textbf{The case} $x_1\ne0,x_2\ne0, x_3=0$.  In this case $x_1N_1+x_2N_2=0$, but $N_1,N_3$ and
$N_2,N_3$ are \emph{linearly independent} in pairs. Rewrite the system for this case as follows.
$$
\begin{array}{l}
(1,1) \quad  \mu_1(\sigma_{13}\mu_3-\mu_1)=0,\\[1ex]
(2,2) \quad  \mu_2(\sigma_{23}\mu_3-\mu_2)=0,\\[1ex]
(3,3) \quad  \equiv 0
\end{array}
$$
$$
\begin{array}{ll}
(1,2)& \mu_1\mu_3(1-\sigma_{23})-\mu_2\mu_3(1-\sigma_{13})+\mu_1(\sigma_{23}\mu_3-\mu_2)x_1^2-\\
&\hspace{7cm} \mu_2(\sigma_{13}\mu_3-\mu_1)x_2^2=0,\\[1ex]
(2,3)&
x_1^2\mu_2(\sigma_{23}\mu_2-\mu_3)+\mu_1\mu_2(1-\sigma_{31}-\mu_1\mu_3(1-\sigma_{21})+\\
&\hspace{7cm} \mu_2(\sigma_{13}\mu_1-\mu_3)x_2^2=0,\\[1ex]
(3,1)&
-x_2^2\mu_1(\sigma_{13}\mu_1-\mu_3+\mu_2\mu_3(1-\sigma_{12})-\mu_1\mu_2(1-\sigma_{32})-\\
&\hspace{7cm}\mu_1(\sigma_{23}\mu_2-\mu_3)x_1^2=0.
\end{array}
$$
Set $\mu_1=\mu_2=0$. Then the system is fulfilled for arbitrary $\mu_3$. The case $\mu_3=0$ is
already considered. The case $\mu_3\ne0$ gives the $S\cap\{e_1,e_2\}_R$ it 3-{rd} row of the Table
1. \\
Set $\mu_1=0, \mu_2\ne0$. Then
$\sigma_{12}=\mu_2\mu_3,\sigma_{13}=\mu_2\mu_3,\sigma_{23}=-\mu_2\mu_3$. The equation (2,2) yields
$-\mu_2^2(\mu_3^2+1)=0$. The contradiction.\\
Set $\mu_1\ne0, \mu_2=0$. Then
$\sigma_{12}=\mu_1\mu_3,\sigma_{13}=-\mu_1\mu_3,\sigma_{23}=\mu_1\mu_3$. The equation (1,1) yields
$-\mu_1^2(\mu_3^2+1)=0$. The contradiction.\\
Set $\mu_1\ne0, \mu_2\ne0$. Then $\mu_1=\sigma_{13}\mu_3, \mu_2=\sigma_{23}\mu_3$ and the
substitution into (1,2) yields
$$
\mu_3^3(\mu_2-\mu_1)=0.
$$
The case $\mu_3=0$ contradicts $\mu_1\ne0, \mu_2\ne0$, as one can see from (1,1) and (2,2). Thus,
set $\mu_1=\mu_2=\mu\ne0$. Then $\sigma_{13}=\sigma_{23}=\mu^2$ and from (1,1) and (2,2) we
conclude
\begin{equation}\label{eqno:14}
\mu\mu_3-1=0.
\end{equation}
In this case we have
\begin{equation}\label{eqno:15}
\sigma_{12}=2-\mu^2,\quad \sigma_{13}=\mu^2,\quad \sigma_{23}=\mu^2.
\end{equation}
The substitution of \eqref{eqno:14} and \eqref{eqno:15} into the system yields the identity. Since
$\mu\mu_3=1$ in our consideration means $\rho_1=\rho_2=2$, we get the 8-{th} row of the Table 1.

\textbf{The case} $x_1\ne0,x_2=0, x_3\ne0$, after similar computations,  resulting
$S\cap\{e_1,e_3\}_R$ in the {3-rd row } and the {9-th row of the Table 1.}

\textbf{The case} $x_1=0,x_2\ne0, x_3\ne0$ resulting $S\cap\{e_2,e_3\}_R$ in the 4-rd row  and the
10-th row of the Table 1.

\textbf{The case} $x_1=1,x_2=0, x_3=0$.  In this case $N_1=0$ and the equations (1,1), (2,2), (3,3)
and (2,3) are fulfilled regardless the geometry of the group. The equations (1,2) and (1,3) take
the forms
$$
\begin{array}{l}
(1,2)\quad \mu_1\mu_3(1-\sigma_{23})-\mu_2\mu_3(1-\sigma_{13})+\mu_1(\sigma_{23}\mu_3-\mu_2)=0\\[1ex]
(1,3)\quad \mu_2\mu_3(1-\sigma_{12})-\mu_1\mu_2(1-\sigma_{23})-\mu_1(\sigma_{23}\mu_2-\mu_3)=0
\end{array}
$$
After simplifications, we get
$$
\begin{array}{l}
(1,2)\quad \sigma_{13}(\mu_2\mu_3-1)=0,\\[1ex]
(1,3)\quad \sigma_{12}(\mu_2\mu_3-1)=0.
\end{array}
$$
The case  $\mu_2\mu_3=1$ means $\rho_1=2$ and we have the 5-th row of the Table 1. Consider the
case $\sigma_{12}=0, \sigma_{13}=0$ which is equivalent to the system
$$\left\{
\begin{array}{l}
 \mu_2\mu_3=0,\\[1ex]
 \mu_1(\mu_2-\mu_3)=0.
\end{array}
\right.
$$
We have 4 possible solutions:
$$(i)\ \mu_1=0,\mu_2=0,\mu_3=0; \quad (ii)\ \mu_1=0,\mu_2=0,\mu_3\ne0;\quad
$$
$$
(iii)\ \mu_1=0,\mu_2\ne0,\mu_3=0;\quad (iv)\ \mu_1\ne0,\mu_2=0,\mu_3=0.
$$
The case (i) is already included into the 1-st row of the Table 1, the case (ii) is already
included into $S\cap\{e_1,e_2\}_R$ case in the 4-st of the Table 1, the case (iii) is already
included into $S\cap\{e_1,e_3\}_R$ case in the 3-rd row of the Table 1. The case (iv) is a new one
and yields $e_1$ field in the 2-nd row of the Table 1.

 \textbf{The case} $x_1=0,x_2=1, x_3=0$ yields $e_2$ into the 3-rd and 6-th rows of the Table 1.

 \textbf{The case} $x_1=0,x_2=0, x_3=1$ yields $e_3$ into the 4-th  and 7-th rows of the Table 1.

The proof is complete.
\end{proof}
Now we specify the result of the Theorem \ref{Thrm:2.1} to each of the unimodular groups.
\begin{theorem}\label{Thrm:2.2}
Let $G$ be a three-dimensional unimodular Lie group with the left-invariant metric and let
$\{e_i,\, i=1,2,3\}$ be an orthonormal basis for the Lie algebra satisfying \eqref{eqno:5}.
Moreover, assume that $\lambda_1\ge\lambda_2\ge\lambda_3$. Then the left-invariant unit vector
fields of G are given as follows:
\begin{center}
\begin{tabular}{|p{0.8in}|c|p{2.3cm}|}
  \hline
  \vphantom{\Bigg|}G & Conditions on $\lambda_1,\lambda_2, \lambda_3$ &  The sets of left-invariant totally geodesic unit vector fields \\
  \hline
  SU(2) & \vphantom{\Big|}$\lambda_1=\lambda_2=\lambda_3=2$ & S \\
   \cline{2-3}& \vphantom{\Big|}$\lambda_1=\lambda_2=\lambda>\lambda_3=2$ & $\pm e_3$\\
  \cline{2-3} & \vphantom{\Bigg|}$\lambda_1=\lambda_2=\lambda>2>\lambda_3=\lambda-\sqrt{\lambda^2-4}$ & $S\cap\{e_1,e_2\}_R$ \\
   \cline{2-3}& \vphantom{\Big|} $\lambda_1=2>\lambda_2=\lambda_3=\lambda>0$ & $\pm e_1$ \\
  \cline{2-3} & \vphantom{\Bigg|} $\lambda_1=\lambda+\sqrt{\lambda^2-4}>\lambda=\lambda_2=\lambda_3>2$ & $S\cap\{e_2,e_3\}_R$ \\
  \cline{2-3} & \vphantom{\Big|} $\lambda_1>\lambda_2>\lambda_3>0$, \ \ $\lambda_m^2-(\lambda_i-\lambda_k)^2=4$ & $\pm e_m$ (i,k,m=1,2,3) \\
   \hline
  SL(2,R) & \vphantom{\Big|} $\lambda_3^2-(\lambda_1-\lambda_2)^2=4$ & $\pm e_3$\\
  \cline{2-3}& \vphantom{\Big|}$ \lambda_1^2-(\lambda_2-\lambda_3)^2=4$ & $\pm e_1$\\
  \hline
   E(2)  & \vphantom{\Big|}$\lambda_1=\lambda_2>0, \ \ \lambda_3=0$ & $\pm e_3$, $S\cap\{e_1,e_2\}_R$ \\
 \cline{2-3} & \vphantom{\Big|}$\lambda_1^2-\lambda_2^2=4$, $\lambda_1>\lambda_2>0$,\ \ $\lambda_3=0$  & $\pm e_1$ \\
 \hline
 E(1,1) &\vphantom{\Big|} $\lambda_1^2-\lambda_2^2=-4$, $\lambda_1>0, \lambda_2<0, \ \ \lambda_3=0$ & $\pm e_2$ \\
  \cline{2-3} &\vphantom{\Big|} $\lambda_1^2-\lambda_2^2=4$, $\lambda_1>0, \lambda_2<0, \ \ \lambda_3=0$ & $\pm e_1$  \\
  \hline
 {\rm Heisenberg group} & \vphantom{\Big|}$\lambda_1=2, \ \ \lambda_2=0,\lambda_3=0$ & $\pm e_1$ \\
 \hline
  $R\oplus R\oplus R$ & \vphantom{\Big|}$\lambda_1=\lambda_2=\lambda_3=0$ & S \\
  \hline
\end{tabular}
\end{center}
where $S\cap\{e_i,e_k\}_R$ means the set of unit vectors in a plane, spanned by $e_i$ and $e_k$, in
the tangent space of the group at the unit element.
\end{theorem}

\begin{proof}

\noindent\textbf{The case SU(2)}. In this case $\lambda_1\geq \lambda_2\geq\lambda_3>0$. A simple
calculation yields
$$
\rho_m=2\mu_i\mu_k=\frac12(\lambda_m^2-(\lambda_i-\lambda_k)^2).
$$
Observe that $\rho_m=\rho_k$ if and only if $\lambda_m=\lambda_k$. From the Table 1 we now readout
the cases
\begin{itemize}
\item if $\lambda_1=\lambda_2=\lambda_3=2$, then $\rho_1=\rho_2=\rho_3=2$ and each left-invariant unit vector
field is totally geodesic one.

\item if $\lambda_1=\lambda_2=\lambda>\lambda_3=2$, then
$\rho_1=\rho_2=\frac12(\lambda^2-(\lambda-2)^2)=2(\lambda-1)>2$, \ $\rho_3=2$ and we have $\pm e_3$
as a unique totally geodesic left-invariant unit vector field.

\item if $\lambda_1=\lambda_2=\lambda>2>\lambda_3$, then
$\rho_1=\rho_2=\frac12(\lambda^2-(\lambda-\lambda_3)^2)=\frac12(2\lambda\lambda_3-\lambda_3^2)$, \
$\rho_3=\frac12 \lambda_3^2<2$. Equalizing
$$
\frac12(2\lambda\lambda_3-\lambda_3^2)=2
$$
we have $\lambda_3=\lambda\pm \sqrt{\lambda^2-4}$. Since $\lambda_3<\lambda$, the appropriate
solution is $\lambda_3=\lambda-\sqrt{\lambda^2-4}$. In this case the set of totally geodesic
left-invariant unit vector fields is $S\cap\{e_1,e_2\}_R$.

\item if $\lambda_1=2>\lambda_2=\lambda_3=\lambda>0$, then $\rho_1=2$,
$\rho_2=\rho_3=\frac12(\lambda^2-(\lambda-2)^2)=2(\lambda-1)<2$ and we have a unique left-invariant
totally geodesic unit vector field $\pm e_1$.

\item if $\lambda_1>2>\lambda_2=\lambda_3=\lambda>0$, then $\rho_1=\frac12\lambda_1^2>2$ and
$\rho_2=\rho_3=\frac12(\lambda^2-(\lambda-\lambda_1)^2)=\frac12(2\lambda\lambda_1-\lambda_1^2)$.
Equalizing
$$
\frac12(2\lambda\lambda_1-\lambda_1^2)=2,
$$
we find $\lambda_1=\lambda\pm \sqrt{\lambda_2-4}$. Since $\lambda_1>\lambda$, the appropriate
solution is $\lambda_1=\lambda+\sqrt{\lambda^2-4}$. In this case the set of totally geodesic
left-invariant unit vector fields is $S\cap\{e_2,e_3\}_R$.

\item if $\lambda_1>\lambda_2>\lambda_3>0$, then $\rho_1$, $\rho_2$ and $\rho_3$ are all different. In this
case, if $$\lambda_m^2-(\lambda_i-\lambda_k)^2=4$$ for $m\ne i\ne k$, then the corresponding Ricci
curvature $\rho_m=2$ and we have $\pm e_m$ as a unique left-invariant totally geodesic unit vector
field.
\end{itemize}

\noindent\textbf{The case SL(2,R)}. In this case $\lambda_1\geq\lambda_2>0$, $\lambda_3<0$ and the
Ricci principal curvatures are
$$
\rho_m=2\mu_i\mu_k=\frac12(\lambda_m^2-(\lambda_i-\lambda_k)^2).
$$
\begin{itemize}
\item if $\lambda_1=\lambda_2=\lambda>0$, then
$\rho_1=\rho_2=\frac12(\lambda^2-(\lambda-\lambda_3)^2)=\frac12(2\lambda\lambda_3-\lambda_3^2)$, \
$\rho_3=\frac12 \lambda_3^2$. Equalizing
$$
\frac12(2\lambda\lambda_3-\lambda_3^2)=2
$$
we have $\lambda_3=\lambda\pm \sqrt{\lambda^2-4}$. Since $\lambda_3<0$, we have no appropriate
solutions. Therefore, equalizing $\rho_3=2$, we have a unique case $\lambda_3=-2$ and the vector
field $\pm e_3$.

\item if $\lambda_1>\lambda_2>0$, then $\rho_1$, $\rho_2$ and $\rho_3$ are all different. In this case, consider
separately the condition
$$
\lambda_m^2-(\lambda_i-\lambda_k)^2=4
$$
for each $m, i, k$.

For $m=3$ we have
$$
\lambda_3^2-(\lambda_1-\lambda_2)^2=4.
$$
If $\lambda_1, \lambda_2$ and $\lambda_3$ satisfy this equation, then $\pm e_3$ is totally
geodesic. Remark, that this case contains the case $\lambda_1=\lambda_2$.

For $m=2$ we have
$$
\lambda_2^2-(\lambda_1-\lambda_3)^2=4.
$$
Since $\lambda_3<0$ we have $\lambda_1-\lambda_3>\lambda_1$. Therefore,
$\lambda_2^2-(\lambda_1-\lambda_3)^2<0$. This contradiction shows that $\pm e_2$ is never totally
geodesic.

For $m=1$ we have
$$
\lambda_1^2-(\lambda_2-\lambda_3)^2=4.
$$
Since $\lambda_1>\lambda_2$ we have $\pm e_1$ totally geodesic  for all solutions of the equation
above. Remark, that the solution necessarily satisfy $\lambda_1-\lambda_2>-\lambda_3$.
\end{itemize}

\noindent\textbf{The case E(2)}. In this case $\lambda_1\geq\lambda_2>0$, $\lambda_3=0$ and the
Ricci principal curvatures are
$$
\rho_1=\frac12(\lambda_1^2-\lambda_2^2), \quad
\rho_2=-\rho_1=\frac12(\lambda_2^2-\lambda_1^2),\quad \rho_3=-\frac12(\lambda_1-\lambda_2)^2.
$$
\begin{itemize}
\item if $\lambda_1=\lambda_2=\lambda>0$ then $\rho_1=\rho_2=\rho_3=0$ and the group is flat. Make an auxiliary
calculations:
$$
\mu_1=\frac12(-\lambda_1+\lambda_2+\lambda_3)=0,\quad
\mu_2=\frac12(\lambda_1-\lambda_2+\lambda_3)=0,
$$
$$
\mu_3=\frac12(\lambda_1+\lambda_2-\lambda_3)=\lambda>0.
$$
From the Table 1 we find $\pm e_3, S\cap\{e_1,e_2\}_R$.

\item if $\lambda_1>\lambda_2$, then we have one more condition $\rho_1=2$,i.e.
$$ \lambda_1^2-\lambda_2^2=4
$$
which yields $\pm e_1$ as a totally geodesic field.

\end{itemize}

\noindent\textbf{The case E(1,1)}. In this case $\lambda_1>0, \lambda_2<0$, $\lambda_3=0$ and the
Ricci principal curvatures are
$$
\rho_1=\frac12(\lambda_1^2-\lambda_2^2), \quad
\rho_2=-\rho_1=\frac12(\lambda_2^2-\lambda_1^2),\quad \rho_3=-\frac12(\lambda_1-\lambda_2)^2.
$$
In this case $\rho_3<0$, $\rho_1\ne\rho_2$ and we have only two possible cases: either $\rho_1=2$
or $\rho_2=2$.
\begin{itemize}
\item $\rho_1=2$. In this case $\lambda_1$ and $\lambda_2$ should satisfy
$$
\lambda_1^2-\lambda_2^2=2,
$$
which yields $\pm e_1$ as the totally geodesic field.

\item $\rho_2=2$. In this case $\lambda_1$ and $\lambda_2$ should satisfy
$$
 \lambda_2^2-\lambda_1^2=4
$$
which yields $\pm e_2$ as a totally geodesic field.
\end{itemize}

\noindent\textbf{The case of Heisenberg group}.  In this case $\lambda_1>0, \lambda_2=\lambda_3=0$,
and the Ricci principal curvatures are
$$
\rho_1=\frac12\lambda_1^2, \quad \rho_2=\rho_3=-\frac12\lambda_1^2.
$$
In this case $\rho_2<0$, $\rho_3<0$ and we have only one possible case $\rho_1=2$ for
$\lambda_1=2$, which yields $\pm e_1$ as the totally geodesic field.\vspace{2ex}

\noindent\textbf{The case $\bf{R\oplus R\oplus R}$}. Here $\lambda_1=\lambda_2=\lambda_3$ and
evidently all left-invariant unit vector fields are totally geodesic.\vspace{1ex}

\end{proof}

\subsection{Geometrical characterization of totally geodesic unit vector fields}

Let $M$ be an odd-dimensional smooth manifold. Denote by $\phi, \xi, \eta$ a $(1,1)$ tensor field ,
a vector field and a 1-form on $M$ respectively. A triple $(\phi,\xi,\eta)$ is called an
\emph{almost contact structure} on $M$ if
\begin{equation}\label{eqno:16}
\phi^2 X=-X+\eta(X)\xi,\quad \phi\xi=0,\quad \eta(\xi)=1
\end{equation}
for any vector field $X$ on $M$. The manifold $M$ with the almost contact structure is called an
\emph{almost contact} manifold.

If $M$ is endowed with a Riemannian metric $g(\cdot,\cdot)=\la\cdot,\cdot\ra$ such that
\begin{equation}\label{eqno:17}
\la\phi X,\phi Y\ra=\la X, Y\ra-\eta(X)\eta(Y),\quad \eta(X)=\la\xi, X\ra
\end{equation}
for all vector fields $X$ and $Y$ on $M$, then a quadruple $(\phi,\xi,\eta,g)$ is called an
\emph{almost contact \textbf{metric} structure} and the manifold is called an \emph{almost contact
\textbf{metric} manifold.} The first of the conditions above is called a \emph{compatibility
condition} for $\phi$ and $g$.

If 2-form $d\eta$, given by
$$
d\eta(X,Y)=\frac12\big(X\eta(Y)-Y\eta(X)-\eta([X,Y])\big),
$$
satisfies
\begin{equation}\label{eqno:18}
d\eta(X,Y)=\la X,\phi Y\ra,
\end{equation}
then the structure $(\phi,\xi,\eta,g)$ is called \textbf{contact metric structure} and the manifold
with a contact metric structure is called by a  \textbf{contact metric manifold}. A contact metric
manifold is called $K$-\emph{contact}, if $\xi$ is a Killing vector field.

The Nijenhuis torsion of tensor field $T$ of type $(1,1)$ is given by
$$
[T,T](X,Y)=T^2[X,Y]+[TX,TY]-T[TX,Y]-T[X,TY]
$$
and defines a $(1,2)$ tensor field on $M$. An almost contact structure  $(\phi,\xi,\eta)$ is called
\emph{normal}, if
\begin{equation}\label{eqno:19}
[\phi,\phi](X,Y)+2d\eta(X,Y)\,\xi=0. \
\end{equation}
Finally, a contact metric structure $(\phi,\xi,\eta,g)$ is called \emph{Sasakian}, if it is normal.
A manifold with Sasakian structure is called \emph{Sasakian manifold}. In Sasakian manifold
necessarily $\phi=A_\xi$ and $ \eta=\la\xi,\cdot \ra$. The unit vector field $\xi$ is called a
\emph{characteristic vector field} of the Sasakian structure and is a Killing one. This vector
field is always totally geodesic \cite{Ym3}.

In tree-dimensional case we have
\begin{theorem}\label{Thrm:2.3}\cite{Ym3}
Let $\xi$ be a unit Killing  vector field on 3-dimensional Riemannian manifold $M^3$. If $\xi(M^3)$
is totally geodesic in $T_1M^3$ then  either
$$
\big(\phi=A_\xi,\ \xi ,\ \eta=\la\xi,\,\cdot\,\ra\big)
$$
is a Sasakian structure on $M^3$ or $M^3=M^2\times E^1$ metrically and $\xi$ is the unit vector
field of Euclidean factor.
\end{theorem}
 Define the structure
\begin{equation}\label{eqno:20}
\big(\phi=A_\xi,\ \xi ,\ \eta=\la\xi,\,\cdot\,\ra\big),
\end{equation}
where the (1,1) tensor field is given by \eqref{W}. Now we can give a geometrical description of
totally geodesic unit vector fields.
\begin{proposition}
Let $\xi$ be a left invariant totally geodesic unit vector field on $SU(2)$ with the left invariant
metric $g$ and let $\{e_i,\, i=1,2,3\}$ be an orthonormal basis for the Lie algebra satisfying
\eqref{eqno:5}. Assume in addition that $\lambda_1\ge\lambda_2\ge\lambda_3$. Then
$$
\big(\phi=A_\xi,\ \xi ,\ \eta=\la\xi,\,\cdot\,\ra\big)
$$
is the almost contact structure on $SU(2)$. Moreover,
\begin{itemize}
\item if $\lambda_1=\lambda_2=\lambda_3=2$ or $\lambda_1=\lambda_2>\lambda_3=2$ or $\lambda_1=2>\lambda_2=\lambda_3$,
      then the structure is Sasakian;
\item if $\lambda_2=\lambda_2=\lambda>2>\lambda_3=\lambda-\sqrt{\lambda^2-4}$ or
$\lambda_1=\lambda+\sqrt{\lambda^2-4}>\lambda=\lambda_2=\lambda_3>2$ , then the structure is
neither normal nor metric;
\item if $\lambda_1>\lambda_2>\lambda_3$, then the structure is normal only for
$$\xi=e_1,\quad \lambda_1=\lambda_2+\frac{1}{\lambda_2},\quad\lambda_3=\frac{1}{\lambda_2},\quad \lambda_2>1$$
\end{itemize}
\end{proposition}
\begin{proof}
Consider the cases from Theorem \ref{Thrm:2.2}.

$\bullet$\ In the case of  $\lambda_1=\lambda_2=\lambda_3=2$ we have $\mu_1=\mu_2=\mu_3=1$ and
hence
$$
A_\xi=\left(
        \begin{array}{ccc}
          0 &  -x_3 &  x_2 \\
          x_3 & 0 &  -x_1 \\
          -x_2 &  x_1 & 0 \\
        \end{array}
      \right).
$$
Therefore, the field $\xi$ is the Killing one. By Theorem \ref{Thrm:2.3}, the structure
\eqref{eqno:20} is Sasakian.

In the case of $\lambda_1=\lambda_2=\lambda>\lambda_3=2$  we have $
\mu_1=1,\mu_2=1,\mu_3=\lambda-1$ and  $\xi=\pm e_3$. For $\xi=+\, e_3$ we find
$$
A_\xi=\left(
        \begin{array}{ccc}
          0 &1 & 0 \\
          -1 & 0 & 0 \\
          0 & 0 & 0 \\
        \end{array}
      \right)
$$
and see that again $\xi$ is  the Killing unit vector field. Therefore, the structure
\eqref{eqno:20} is Sasakian.

In the case of $\lambda_1=2>\lambda_2=\lambda_3=\lambda>0$  we have $
\mu_1=-1+\lambda,\mu_2=1,\mu_3=1$ and  $\xi=\pm e_1$. For $\xi=+\, e_1$ we find
$$
A_\xi=\left(
        \begin{array}{ccc}
          0 &0 & 0 \\
          0 & 0 & -1 \\
          0 & 1 & 0 \\
        \end{array}
      \right)
$$
and see that again $\xi$ is  the Killing unit vector field. Therefore, the structure
\eqref{eqno:20} is Sasakian.

$\bullet$\ Consider the case $\lambda_1=\lambda_2=\lambda>2>\lambda_3=\lambda-\sqrt{\lambda^2-4}$
and $\xi=x_1e_1+x_2e_2$. We have
$$
\mu_1=\frac12(\lambda-\sqrt{\lambda^2-4}), \quad \mu_2=\frac12 (\lambda-\sqrt{\lambda^2-4}),\quad
\mu_3=\frac12(\lambda+\sqrt{\lambda^2-4}).
$$
Set for brevity $\theta =\frac12(\lambda-\sqrt{\lambda^2-4})$ and $\bar
\theta=\frac12(\lambda+\sqrt{\lambda^2-4})$. Then
$$
\mu_1=\theta, \quad \mu_2=\theta ,\quad \mu_3=\bar\theta, \quad \theta\bar\theta=1 \quad (\theta\ne
1,\bar\theta\ne 1)
$$
and for this case we have
$$
A_\xi=\left(
        \begin{array}{ccc}
          0 & 0 & \bar\theta x_2 \\
          0 & 0 & -\bar\theta x_1 \\
          -\theta x_2 & \theta x_1 & 0 \\
        \end{array}
      \right).
$$
Since $\theta\ne\bar\theta$, the field $\xi$ is never Killing one but geodesic. Indeed,
$$
A_\xi\xi=\left(
        \begin{array}{ccc}
          0 & 0 & \bar\theta x_2 \\
          0 & 0 & -\bar\theta x_1 \\
          -\theta x_2 & \theta x_1 & 0 \\
        \end{array}
      \right)
      \left(
        \begin{array}{c}
          x_1 \\
          x_2 \\
          0 \\
        \end{array}
      \right)=
            \left(
        \begin{array}{c}
          0 \\
          0 \\
          \theta(-x_2x_1+x_1x_2) \\
        \end{array}
      \right)=0.
$$
The structure \eqref{eqno:20} is an almost contact one on $SU(2)$. Indeed,
\begin{multline*}
\phi^2=\left(
        \begin{array}{ccc}
          0 & 0 & \bar\theta x_2 \\
          0 & 0 & -\bar\theta x_1 \\
          -\theta x_2 & \theta x_1 & 0 \\
        \end{array}
      \right)
      \left(
        \begin{array}{ccc}
          0 & 0 & \bar\theta x_2 \\
          0 & 0 & -\bar\theta x_1 \\
          -\theta x_2 & \theta x_1 & 0 \\
        \end{array}
      \right)=
\\[1ex]
    \left(
        \begin{array}{ccc}
          -x_2^2 & x_1x_2 & 0 \\
          x_1x_2 & -x_1^2 & 0 \\
          0 & 0 & -1 \\
        \end{array}
      \right).
\end{multline*}
Then
$$
\phi^2Z=\left(
        \begin{array}{ccc}
          -1+x_1^2 & x_1x_2 & 0 \\
          x_1x_2 & -1+x_2^2 & 0 \\
          0 & 0 & -1 \\
        \end{array}
      \right)
      \left(
        \begin{array}{c}
          z_1 \\
          z_2 \\
          z_3 \\
        \end{array}
      \right)=-Z+\la \xi, Z\ra \xi.
$$
This structure is not metric one.  For the  compatibility condition \eqref{eqno:17} we have
$$
\phi Z=\left(
         \begin{array}{c}
           \bar\theta x_2z_3  \\
           -\bar\theta x_1z_3 \\
           -\theta x_2z_1+\theta x_1z_2 \\
         \end{array}
       \right),\quad
\phi W =\left(
         \begin{array}{c}
           \bar\theta x_2w_3  \\
           -\bar\theta x_1w_3 \\
           -\theta x_2w_1+\theta x_1w_2 \\
         \end{array}
       \right)
$$
and hence
\begin{multline*}
\la \phi Z,\phi
W\ra=\bar\theta^2z_3w_3+\theta^2(x_2^2z_1w_1+x_1^2z_2w_2-x_1x_2z_1w_2-x_1x_2z_2w_1)=\\
\theta^2(z_1w_1+z_2w_2)+\bar\theta^2x_3w_3-\la\xi,Z\ra\la\xi,W\ra\ne \la Z,
W\ra-\la\xi,Z\ra\la\xi,W\ra.
\end{multline*}
This structure is not normal one. To prove this, check the normality condition \eqref{eqno:16}.
Find the Nijenhuis torsion of $\phi$ on $e_1,e_2$. We have
$$
\begin{array}{l}
\phi e_1=-\theta x_2 e_3,\quad \phi e_2=\theta x_1 e_3,\quad \phi e_3=\bar\theta x_2
e_1-\bar\theta x_1 e_2,\\[1ex]
[e_1,e_2]=2\theta e_3,\quad [e_1,e_3]=-(\theta+\bar\theta) e_2,\quad
[e_2,e_3]=(\theta+\bar\theta)e_1, \\[1ex]
\phi^2[e_1,e_2]=-2\theta e_3, \quad [\phi e_1,\phi e_2]=0, \\[1ex]
\phi [\phi e_1,e_2]=-\theta^2x_2^2(\theta+\bar\theta) e_3=-\theta (\theta^2+1)x_2^2 e_3,
\\[1ex]
\phi [e_1,\phi e_2]=-\theta^2(\theta+\bar\theta)x_1^2e_3=-\theta (\theta^2+1)x_1^2 e_3
\end{array}
$$
and thus,
$$
[\phi,\phi](e_1,e_2)=\theta(\theta^2-1)e_3\ne 2d\eta(e_1,e_2) \xi.
$$
In a similar way we can analyze the case
$\lambda_1=\lambda+\sqrt{\lambda^2-4}>\lambda=\lambda_2=\lambda_3>2$ with the same result.

$\bullet$\ Consider the case $\lambda_1>\lambda_2>\lambda_3$, $\xi=\pm e_i$.  We have
$$
\mu_1=\frac12(-\lambda_1+\lambda_2+\lambda_3), \quad
\mu_2=\frac12(\lambda_1-\lambda_2+\lambda_3),\quad \mu_3=\frac12(\lambda_1+\lambda_2-\lambda_3)
$$
Set $\xi=e_1$. The condition  $\lambda_1^2-(\lambda_2-\lambda_3)^2=4$ means that $\mu_2\mu_3=1$.
The matrix $A_\xi$ takes the form
$$
A_\xi=\left(
        \begin{array}{ccc}
          0 & 0 & 0 \\
          0 & 0 & -\mu_3 \\
         0& \mu_2  & 0 \\
        \end{array}
      \right).
$$
Since $\mu_2\ne \mu_3$, the field $\xi$ is not a Killing one, but geodesic. The structure
\eqref{eqno:20} is almost contact one. Indeed,
$$
\phi^2=\left(
        \begin{array}{ccc}
          0 & 0 & 0 \\
          0 & -\mu_2\mu_3 &0 \\
         0& & -\mu_3\mu_2 \\
        \end{array}
      \right)=\left(
        \begin{array}{ccc}
          0 & 0 & 0 \\
          0 & -1 & 0 \\
         0& 0  & -1 \\
        \end{array}
      \right)
$$
and hence
$$
\phi^2 Z=-Z+\la \xi,Z\ra \xi.
$$
The structure is \emph{normal} if and only if
\begin{equation}\label{eqno:21}
\lambda_1=\lambda_2+\frac{1}{\lambda_2},\quad \lambda_3=\frac{1}{\lambda_2},\quad  \lambda_2>1.
\end{equation}
Indeed, remark that
$$
\phi e_1=0, \quad \phi e_2=\mu_2 e_3,\quad \phi e_3=-\mu_3 e_2.
$$
Now set $Z=e_1, W=e_2$. Then we have
$$
\begin{array}{l}
\phi^2[e_1,e_2]=\lambda_3\phi^2e_3=-\lambda_3e_3, \\[1ex]
[\phi e_1,\phi e_2]=0, \\[1ex]
\phi[\phi e_1,e_2]=0, \\[1ex]
\phi[e_1,\phi e_2]=\mu_2\phi[e_1,e_3]=-\mu_2\lambda_2\phi e_2=-\mu_2^2\lambda_2e_3, \\[1ex]
d\eta (e_1,e_2)=\frac12(\la e_1,\phi e_2\ra-\la\phi e_1,e_2\ra)=0.
\end{array}
$$
Therefore, the first necessary condition of normality is $ \lambda_3=\mu_2^2\lambda_2 $. If we
remark that $\mu_2\mu_3=1$, we can rewrite this condition as
\begin{equation}\label{eqno:22}
\lambda_3\mu_3=\lambda_2\mu_2.
\end{equation}
Set $Z=e_1, W=e_3$. Then we have
$$
\begin{array}{l}
\phi^2[e_1,e_3]=-\lambda_2\phi^2e_2=\lambda_2e_2, \\[1ex]
[\phi e_1,\phi e_3]=0, \\[1ex]
\phi[\phi e_1,e_3]=0, \\[1ex]
\phi[e_1,\phi e_3]=-\mu_3\phi[e_1,e_2]=-\mu_3\lambda_3\phi e_3=\mu_3^2\lambda_3e_2, \\[1ex]
d\eta (e_1,e_3)=\frac12(\la e_1,\phi e_3\ra-\la\phi e_1,e_3\ra)=0.
\end{array}
$$
Therefore, the second necessary condition of normality is $ \lambda_2=\mu_3^2\lambda_3 $. which is
equivalent to \eqref{eqno:22}.

Finally, set $Z=e_2, W=e_3$. Then we have
$$
\begin{array}{l}
\phi^2[e_2,e_3]=\lambda_1\phi^2e_1=0, \\[1ex]
[\phi e_2,\phi e_3]=-\mu_2\mu_3[e_3, e_2]=\lambda_1e_1, \\[1ex]
\phi[\phi e_2,e_3]=0, \\[1ex]
\phi[e_2,\phi e_3]=0, \\[1ex]
d\eta (e_2,e_3)=\frac12(\la e_2,\phi e_3\ra-\la\phi
e_2,e_3\ra)=-\frac12(\mu_3+\mu_2)=-\frac12\lambda_1.
\end{array}
$$
These data satisfy \eqref{eqno:19}. Expand the equation \eqref{eqno:22}, namely
$$
\lambda_3(\lambda_1+\lambda_2-\lambda_3)=\lambda_2(\lambda_1-\lambda_2+\lambda_3)
$$
and perform rearrangements as follows:
$$
\lambda_1(\lambda_3-\lambda_2)+\lambda_3(\lambda_2-\lambda_3)=\lambda_2(-\lambda_2+\lambda_3).
$$
Since  $\lambda_2\ne \lambda_3$, we get
$$
\lambda_1=\lambda_2+\lambda_3.
$$
Then
$$
\mu_1=0,\quad \mu_2=\frac12(\lambda_1-\lambda_2+\lambda_3)=\lambda_3,\quad
\mu_3=\frac12(\lambda_1+\lambda_2-\lambda_3)=\lambda_2
$$
and, from the condition $\mu_2\mu_3=1$,  we find
$$
\lambda_2\lambda_3=1.
$$
Since $\lambda_1>\lambda_2>\lambda_3$, we get \eqref{eqno:21}.

The structure \emph{is not metric}, since
$$
\la \phi Z,\phi W\ra=\mu_3^2z_3w_3+\mu_2^2z_2w_2\ne\la Z,
W\ra-\la\xi,Z\ra\la\xi,W\ra=z_2w_2+z_3w_3.
$$

Setting $\xi=e_2$, we get the normality condition of the form $\lambda_2=\lambda_1+\lambda_3$ which
contradicts the condition $\lambda_1>\lambda_2>\lambda_3$. The structure is not metric.

Setting $\xi=e_3$, we get the normality condition of the form $\lambda_3=\lambda_1+\lambda_2$ which
contradicts again the condition $\lambda_1>\lambda_2>\lambda_3$. The structure is not metric.
\end{proof}

\begin{proposition}
Let $\xi$ be a left invariant totally geodesic unit vector field on $SL(2,R)$ with the left
invariant metric $g$ and let $\{e_i,\, i=1,2,3\}$ be an orthonormal basis for the Lie algebra
satisfying \eqref{eqno:5}. Assume in addition that $\lambda_1\ge\lambda_2>0,\lambda_3<0$. Then
$$
\big(\phi=A_\xi,\ \xi ,\ \eta=\la\xi,\,\cdot\,\ra\big)
$$
is the almost  contact structure on $SL(2,R)$, where $\la\cdot\,,\,\cdot\ra$ is the scalar  product
with respect to $g$. Moreover, if
\begin{itemize}
\item $\lambda_1=\lambda_2, \quad\lambda_3=-2$, then the structure is Sasakian;
\item $\lambda_3=-\sqrt{4+(\lambda_1-\lambda_2)^2}<-2$ or
$\lambda_1=\sqrt{4+(\lambda_2-\lambda_3)^2}$ , then the structure is neither
 normal nor metric.
\end{itemize}
\end{proposition}
\begin{proof}
Consider the case of $\lambda_3=-\sqrt{4+(\lambda_1-\lambda_2)^2}\leq -2$ and $\xi=e_3$. We have
$$
\phi=A_\xi=\left(
        \begin{array}{ccc}
          0 & -\mu_2 & 0 \\
          \mu_1& 0 & 0 \\
          0 & 0 & 0 \\
        \end{array}
      \right)
$$
with $\mu_1\mu_2=1$.

If $\lambda_1=\lambda_2$, then $\lambda_3=-2$ and $\mu_1=\mu_2=1$. Thus the field $\xi=e_3$ is the
Killing one and the structure is Sasakian.

If $\lambda_1>\lambda_2$, then $\lambda_3<-2$. The structure is almost contact, since
$$
\phi^2=\left(
        \begin{array}{ccc}
          -1 & 0 & 0 \\
          0& -1 & 0 \\
          0 & 0 & 0 \\
        \end{array}
      \right).
$$
Similar to the $SU(2)$ case, the structure is not metric and the normality condition  for $\xi=e_3$
takes the form $\lambda_3=\lambda_1+\lambda_2$, which contradicts the sign conditions on
$\lambda_i$.

Consider the case $\lambda_1=\sqrt{4+(\lambda_2-\lambda_3)^2}$ and $\xi=e_1$.  We have
$$
\phi=A_\xi=\left(
             \begin{array}{ccc}
               0 & 0 & 0 \\
               0 & 0 & -\mu_3 \\
               0 & \mu_2 & 0 \\
             \end{array}
           \right)
$$
with $\mu_2\mu_3=1$ \ ($\mu_2\ne1,\ \mu_3\ne 1$). Similar to $SU(2)$ case 3, the normality
conditions take the form $\lambda_2=\mu_3^2\lambda_3$ and $\lambda_3=\mu_2^2\lambda_2$, that
contradicts again the sign conditions on $\lambda_i$.
\end{proof}
\begin{proposition}
Let $\xi$ be a left invariant totally geodesic unit vector field on $E(2)$ with the left invariant
metric $g$ and let $\{e_i,\, i=1,2,3\}$ be an orthonormal basis for the Lie algebra satisfying
\eqref{eqno:5}. Assume in addition that $\lambda_1\ge\lambda_2>0,\lambda_3=0$.

If $\lambda_1=\lambda_2=\lambda>0$, then the group is flat. Moreover,
\begin{itemize}
\item if  $\xi=e_3$, then $\xi$ is a parallel vector field on $E(2)$;
\item if $\xi=x_1e_1+x_2e_2$, then $\xi$ moves along $e_3$ with a constant angle speed~$\lambda$.
\end{itemize}
If $\lambda_1>\lambda_2>0$, then $ \big(\phi=A_\xi,\ \xi ,\ \eta=\la\xi,\,\cdot\,\ra\big) $ is the
almost  contact structure on $E(2)$. This structure is neither metric nor normal.
\end{proposition}
\begin{proof}
Set $\lambda_1=\lambda_2=\lambda>0$. Then $\mu_1=0,\mu_2=0,\mu_3=\lambda$ and for $\xi=e_3$ we have
$A_\xi=0$. This means that $\xi$ is a parallel vector field. If $\xi=x_1e_1+x_2e_2$, then
$$
A_\xi=\left(
        \begin{array}{ccc}
          0 & 0 & \lambda x_2 \\
          0 & 0 & -\lambda x_1 \\
          0 & 0 & 0 \\
        \end{array}
      \right)
$$
Since $A_\xi^2=0$, the structure \eqref{eqno:20} is not almost contact one. The field $\xi$ is not
Killing but geodesic one. Moreover,
$$
\nabla_{e_3}\xi=\lambda(x_1e_2-x_2e_1).
$$
This means that the field $\xi$ moves along $e_3$-geodesics with a constant angle speed $\lambda$.

Set $\lambda_1>\lambda_2>0$ and $\xi=e_1$. Then $\mu_1=\frac12(-\lambda_1+\lambda_2)$,
$\mu_2=\frac12(\lambda_1-\lambda_2)$, $\mu_3=\frac12(\lambda_1+\lambda_2)$ and $\mu_2\mu_3=1$. We
have
$$
A_\xi=\left(
        \begin{array}{ccc}
          0 & 0 & 0 \\
          0 & 0 & -\mu_3 \\
          0 & \mu_2 & 0 \\
        \end{array}
      \right)
$$
The structure \eqref{eqno:20} is an almost contact one. Similar to $SU(2)$ case 3 ( with
$\lambda_3=0$ setting), the normality condition of this structure is
$\lambda_2=\mu_3^2\lambda_3\,(=0)$ which yields a contradiction.
\end{proof}
\begin{proposition}
Let $\xi$ be a left invariant totally geodesic unit vector field on $E(1,1)$ with the left
invariant metric and let $\{e_i,\, i=1,2,3\}$ be an orthonormal basis for the Lie algebra
satisfying \eqref{eqno:5}. Assume in addition that $\lambda_1>0,\lambda_2<0,\lambda_3=0$. Then
$$\big(\phi=A_\xi,\ \xi ,\ \eta=\la\xi,\,\cdot\,\ra\big) $$ is the almost  contact structure on
$E(1,1)$. This structure is neither metric nor normal.
\end{proposition}
\begin{proof}
Consider the case $\lambda_1^2-\lambda_2^2=-4$, which is equivalent to $\mu_1\mu_2=1$, and
$\xi=e_3$. Then
$$
A_\xi=\left(
        \begin{array}{ccc}
          0 & -\mu_2 & 0 \\
          \mu_1 & 0 & 0 \\
          0 & 0 & 0 \\
        \end{array}
      \right)
$$
and the structure is almost contact one. As in previous cases, the structure is neither metric nor
normal. A conclusion is true for the case of $\lambda_1^2-\lambda_2^2=4$ and $\xi=e_1$.
\end{proof}
\begin{proposition}
Let $\xi$ be a left invariant totally geodesic unit vector field on Heisenberg group with the left
invariant metric and let $\{e_i,\, i=1,2,3\}$ be an orthonormal basis for the Lie algebra
satisfying \eqref{eqno:5}. Moreover, assume that $\lambda_1>0,\lambda_2=0,\lambda_3=0$. Then
$$\big(\phi=A_\xi,\ \xi ,\ \eta=\la\xi,\,\cdot\,\ra\big) $$ is the Sasakian structure.
\end{proposition}
\begin{proof}
Indeed, for this case we have $\mu_1=-1, \mu_2=1,\mu_3=1$ and $\xi=e_1$. We have
$$
A_\xi=\left(
        \begin{array}{ccc}
          0 & 0 & 0 \\
          0 & 0 & -1 \\
          0 & 1 & 0 \\
        \end{array}
      \right)
$$
which means that $\xi$ is a Killing vector field and thus the structure is Sasakian.
\end{proof}
\section{Non-unimodular case.}
Choose the orthonormal frame $e_1,e_2, e_3$ as in \eqref{eqno:6}.  Then the Levi-Civita connection
is given by the following table
\begin{equation}\label{eqno:23}
\begin{tabular}{|c|c|c|c|}
  \hline
  \vf $\nabla$ & $e_1$ & $e_2$ & $e_3$ \\[1ex]
  \hline
\vf  $e_1$ & 0 & $\bt e_3$ & $ -\bt e_2$ \\
  \hline
\vf  $ e_2 $ & $-\al e_2$ & $\al e_1$ & $0$ \\
  \hline
\vf  $e_3$ & $ -\dt e_3$ & $0$ & $\dt e_1$ \\
  \hline
\end{tabular}
\end{equation}
For any left-invariant unit vector field $\xi=x_1e_1+x_2e_2+x_3e_3$ we have
$$
\nabla_{e_1}\xi=\bt\, e_1\times\xi, \quad \nabla_{e_2}\xi=-\al\,  e_3\times\xi,\quad
\nabla_{e_3}\xi=\dt \, e_2\times\xi.
$$
Set for brevity
$$
e_1\times\xi=N_1, \quad e_3\times\xi=N_2,\quad e_2\times\xi=N_3,
$$
or in explicit form
\begin{equation}\label{eqno:24}
\begin{tabular}{|c|c|c|}
  \hline
  $N_1$ & $N_2$ & $N_3$ \\
  \hline
  0 & $-x_2$ & $\hphantom{-}x_3$ \\
  $-x_3$ & $\hphantom{-}x_1$ & 0 \\
  $\hphantom{-}x_2$ & $0$ & $-x_1$ \\
  \hline
\end{tabular}.
\end{equation}
 Then
$$
A_\xi e_1=-\bt N_1,\quad A_\xi e_2=\al N_2,\quad A_\xi e_3=-\dt N_3
$$
and the matrix of $A_\xi$ takes the form
\begin{equation}\label{eqno:25}
A_\xi=
\left(%
\begin{array}{ccc}
  0 & -\al x_2 & -\dt x_3 \\[1ex]
  \bt x_3 & \al x_1 & 0 \\[1ex]
  -\bt x_2 & 0 & \dt x_1 \\
\end{array}%
\right)
\end{equation}
A direct computation gives the following result.
\begin{lemma}\label{Lemma:3.1} The derivatives $(\nabla_{e_i}A_\xi)e_k$ of the Weingarten operator $A_\xi$ for the left invariant unit vector field are as
in the following table.\vspace{2ex}

\hspace{-0.5cm}\begin{tabular}{|c||c|c|c|}
  \hline
\vf   & $e_1$ & $e_2$ & $e_3$ \\[1ex]
  \hline\hline
\vf  $e_1$ &$-\bt ^2(x_1e_1-\xi)$ & $\bt \dt  N_3 +\bt \al x_1 e_3$ & $ \bt \al N_2-\bt\dt x_1 e_2$  \\[1ex]
  \hline
\vf  $e_2 $ & $\al^2 N_2 +\bt\al x_3 e_1$ & $ \bt\al N_1-\al^2(x_3e_3-\xi)$  & $\al\dt x_3 e_2$\\[1ex]
  \hline
\vf  $e_3$ & $-\dt^2 N_3-\bt \dt x_2e_1$   & $\al\dt x_2 e_3$& $\bt\dt N_1-\dt^2(x_2e_2-\xi)$  \\[1ex]
  \hline
\end{tabular}
\end{lemma}
\begin{proof}
By definition,
$$
(\nabla_{e_i}A_\xi)e_k=\nabla_{\nabla_{e_i}e_k}\xi-\nabla_{e_i}\nabla_{e_k}\xi.
$$
Using the Table \eqref{eqno:23}, we can easily fill out the table
$$
\begin{tabular}{|c||c|c|c|}
  \hline
  \vphantom{\Big|}$\nabla_{\nabla_{e_i}e_k}\xi$ & $e_1$ & $e_2$ & $e_3$ \\[1ex]
  \hline\hline
\vphantom{\Big|}  $e_1$ &0 & $\bt \dt  N_3$ & $ \bt \al N_2$ \\[1ex]
  \hline
\vphantom{\Big|}  $e_2 $ & $\al^2 N_2$ & $ \bt\al N_1 $ & $0$ \\[1ex]
  \hline
\vphantom{\Big|}  $e_3$ & $-\dt^2 N_3$ & $0$ & $\bt\dt N_1$ \\[1ex]
  \hline
\end{tabular}
$$
and the table
$$
\begin{tabular}{|c||c||c|c|c|}
  \hline
  \vphantom{\Big|}$\nabla$ & $\xi$ &$\nabla_{e_1}\xi$ & $\nabla_{e_2}\xi$ & $\nabla_{e_3}\xi$ \\[1ex]
  \hline  \hline
  \vphantom{\Big|}$e_1$ &$\bt N_1$ &$\bt ^2(x_1e_1-\xi)$ & $-\bt \al x_1 e_3$ & $ \bt\dt x_1 e_2$ \\[1ex]
  \hline
  \vphantom{\Big|}$e_2 $ &$-\al N_2$ & $-\bt\al x_3 e_1$ & $\al^2(x_2e_2-\xi)$ & $-\al\dt x_3 e_2$ \\[1ex]
  \hline
  \vphantom{\Big|}$e_3$ &$\dt N_3$ & $\bt \dt x_2e_1$ & $-\al\dt x_2 e_3$ & $\dt^2(x_3e_3-\xi)$ \\[1ex]
  \hline
\end{tabular}
$$
Now, the result follows immediately.
\end{proof}

By the straightforward application of Codazzi equation and Lemma \ref{Lemma:3.1} we can easily
prove the following.
\begin{lemma} The curvature operator of the non-unimodular group with respect to the chosen
frame takes the form
$$
\begin{array}{l}
R(e_1,e_2)\xi=\al^2N_2 +\bt(\al-\dt)N_3,\\[1ex]
R(e_1,e_3)\xi=-\dt^2N_3 -\bt(\al-\dt)\,N_2\\[1ex]
R(e_2,e_3)\xi=\al\dt N_1
\end{array}
$$
\end{lemma}

Now, everything is prepared for the calculation of the components of total geodesity form
\eqref{eqno:4}.
\begin{lemma}\label{Lemma:3.3}
Let $G$ be non-unimodular Lie group with the basis, satisfying \eqref{eqno:6}. Then the
left-invariant unit vector field $\xi=x_1e_1+x_2e_2+x_3e_3$ is totally geodesic if and only if it
satisfies the following equations:
\begin{multline*}
(1,1)\quad \bt x_1\Big\{\big[\bt [1+\al(\al-\dt)]x_2+\al^3 x_3\big]N_2-\\
\big[\bt [1-\dt(\al-\dt)]x_3-\dt^3x_2\big]N_3\Big\}=0,
\end{multline*}

\begin{multline*}
(2,2)\quad \al \Big\{\Big[\bt[1+\al^2(1-x_3^2)]-[\al+\bt^2(\al-\dt)] x_2x_3\big]\Big]N_1+\\
\al \Big[1+\dt^2\Big]x_1x_3N_3\Big\}=0,
\end{multline*}

\begin{multline*}
(3,3)\quad \dt \Big\{\Big[\bt[1+\dt^2(1-x_2^2)]+[\dt-\bt^2(\al-\dt)]x_2x_3\big]\Big]N_1-\\
\dt \Big[1+\al^2\Big]x_1x_2N_2\Big\}=0,
\end{multline*}

\begin{multline*}
(1,2)\quad  \bt x_1\Big[[\al+\bt^2(\al-\dt)]x_2+\bt\al^2x_3\Big]N_1+ \\
\al\Big[\al[1+\al^2(1-x_3^2)]-\bt [1+\al(\al-\dt)]x_2x_3 \Big]N_2+ \\
\Big[\al\dt
\big[\bt\dt(1-x_1^2)-\dt^2x_2x_3+\bt(\al-\dt)(1-x_3^2)\big]+\bt\al(x_3^2-x_1^2)+\bt\dt\Big]N_3=0,
\end{multline*}

\begin{multline*}
(1,3)\quad  \bt x_1\Big[[\dt-\bt^2(\al-\dt)]x_3-\bt\dt^2x_2\Big]N_1- \\
\Big[\al\dt\big[
 \al\bt(1-x_1^2)+\al^2x_2x_3-\bt(\al-\dt)(1-x_2^2)\big]+\bt\al +\bt\dt(x_2^2-x_1^2) \Big]N_2+ \\
\dt\Big[\bt [-1+\dt(\al-\dt)]x_2x_3 -\dt[1+\dt^2(1-x_2^2)] \Big]N_3=0,
\end{multline*}

\begin{multline*}
(2,3)\quad  \Big[\bt\big[\al\dt(\al+\dt)x_2x_3-\bt(\al-\dt)(\al(1-x_3^2)+\dt(1-x_2^2))\big]+\al\dt(x_2^2-x_3^2) \Big]N_1+ \\
\al\dt\Big[1+\al^2\Big]x_1x_3N_2 -\al\dt\Big[1+\dt^2\Big]x_1x_2N_3=0.
\end{multline*}

\end{lemma}
The proof consists of rather long calculations of the corresponding components $TG_\xi(e_i,e_k)$
for various combinations of $(i,k)$, similar to the calculations in the unimodular case.

The analysis of the Lemma \ref{Lemma:3.3} we split into two subcases.

\begin{theorem}\label{thm:1}
Let $G$ be non-unimodular Lie group with the basis \eqref{eqno:6}. Let $\xi$ be a left invariant
unit vector field which does not belong to the unimodular kernel of the Lie algebra at the origin.
Then $\xi$ is never totally geodesic.
\end{theorem}
\begin{proof}
By the hypothesis, $x_1\ne0$. From \eqref{eqno:24} it follows that $N_2\ne 0$,  $N_3\ne0$ and they
are always linearly independent. Moreover, the vectors $N_1$ and $N_3$ are linearly dependent if
and only if $x_3=0$. If $x_3\ne0$, then the equation $(2,2)$ implies $x_3=0$ and we come to a
contradiction.

Set $x_3=0$. If $x_2\ne0$, then $N_1$ and $N_2$ are linearly independent and $(3,3)$ implies
$\dt=0$. In this case we can rewrite $(1,1)$ as $ \bt^2 x_1 x_2(1+\al^2)N_2=0 $ and we have
$\bt=0$. In this case the equation $(1,2)$ takes the form $\al^2(1+\al^2)N_2=0 $ and we have a
contradiction.

Set $x_3=x_2=0$. In this case $\xi=e_1$, $N_1=0$, $N_2=e_2$ and $N_3=-e_3$. The equation $(1,2)$
takes the form  $\al^2(1+\al^2)N_2=0$. Contradiction.
\end{proof}

\begin{theorem}\label{Thrm:3.2}
Let $G$ be non-unimodular Lie group with the basis as above. Let $\xi$ be a left invariant totally
geodesic unit vector field from the unimodular kernel of the Lie algebra at the origin. Then either
$$\bt=\dt=0 \quad \mbox{and}\quad \xi=\pm \,e_3$$
or
$$
\bt=\theta=\pm\,1,\quad \al\dt=-1\quad \mbox{and}\quad
\pm\,\xi=\theta\,\frac{1}{\sqrt{1+\al^2}}\,e_2+\frac{\al}{\sqrt{1+\al^2}}\,e_3,
$$
\end{theorem}
\begin{proof}
Suppose $\xi=x_2e_2+x_3e_3$. Since $x_1=0$, we have $N_1\ne0$ and $N_1$ is linearly independent
with either $N_2$ or $N_3$.

Suppose $\bt=0$. Then $(2,2)$ implies $-\al^2x_2x_3=0$ and we have the following cases.
\begin{itemize}
\item Case $x_3=0$. Then $N_1=\pm\,e_3$, $N_2=\mp\,e_1$, $N_3=0$
and the equation $(1,2)$ takes the form $\al^2(1+\al^2)N_2=0$. Contradiction.

\item Case $x_2=0$.  Then $N_1=\mp\, e_2$, $N_2=0$, $N_3=\pm\,
e_1$. The equation $(1,3)$ then takes the form $-\dt^2(1+\dt^2)N_3=0$ and we should set $\dt=0$. It
is easy to check that if $\bt=\dt=0$, then all equations are fulfilled. Moreover, the field
$\xi=\pm e_3$ becomes a parallel vector field, since $\nabla \xi=0$.
\end{itemize}

Suppose $\bt\ne 0$, $\dt=0$. Then $(1,3)$ implies $\bt\al N_2=0$ and we have $x_2=0$. In this case
$x_3^2=1$ and $(2,2)$ yields $\al\bt N_1=0$. Contradiction.

Suppose $\bt\ne0, \dt\ne 0$.  In this case a direct analysis of the system becomes too complicated.
Fortunately,  we can apply to this case a different method based on the explicit expression for the
second fundamental form of $\xi(M^n)\subset T_1M^n$ \cite{Ym1}.\vspace{1ex}

\emph{Let $\xi$ be a unit vector field on a Riemannian manifold $M^{n+1}$. The components of second
fundamental form of $\xi(M)\subset T_1M^{n+1}$ can be given by
    $$
\begin{array}{ll}
    \tilde \Omega_{\sigma | i j}= &\frac{1}{2}\Lambda_{\sigma i j}
    \Big\{-\big< (\nabla_{e_i}A_\xi)e_j+  (\nabla_{e_j}A_\xi)e_i, f_\sigma \big>+\\[2ex]
    &\hspace{3cm}\lambda_\sigma\left[ \lambda_j \big< R(e_\sigma, e_i)
    \xi, f_j \big> +  \lambda_i  \big<R(e_\sigma, e_j) \xi, f_i \big>\right] \Big\},
\end{array}
   $$
   where  $\Lambda_{\sigma i j}=[(1+\lambda_\sigma^2)(1+\lambda_i^2)(1+\lambda_j^2)]^{-1/2}$,
   $\lambda_0=0, \lambda_1,\dots, \lambda_n$ are the singular values of the matrix
   $A_\xi$ and $e_0, e_1,\dots, e_n; f_1,\dots, f_n$ are the
   orthonormal frames of   singular vectors  $(i,j=0,1,\dots, n;\,\sigma=1,\dots,n)$.
}\vspace{1ex}

Since $x_1=0$, the matrix \eqref{eqno:25} takes the form
$$
A_\xi=
\left(%
\begin{array}{ccc}
  0 & -\al x_2 & -\dt x_3 \\[1ex]
  \bt x_3 & 0 & 0 \\[1ex]
  -\bt x_2 & 0 & 0 \\
\end{array}%
\right)
$$
Denote by $\te_0,\te_1, \te_2; \tf_1,\tf_2$ the orthonormal singular frames of $A_\xi$.  The matrix
$A_\xi^tA_\xi$ takes the form
\begin{equation}\label{eqno:26}
A_\xi^tA_\xi=
\left(%
\begin{array}{ccc}
  \bt^2 & 0 & 0 \\[1ex]
  0 & \al^2x_2^2 & \al\dt x_2x_3 \\[1ex]
  0 & \al\dt x_2x_3 & \dt^2 x_3^2 \\
\end{array}%
\right).
\end{equation}
The eigenvalues  are $\big[0,\bt^2,\al^2x_2^2+\dt^2x_3^2\big]$. Denote
$m=\sqrt{\vphantom{\big|}\al^2x_2^2+\dt^2x_3^2}$. Then the singular values are
$$
\lambda_0=0,\ \lambda_1=|\bt|,\ \lambda_2=m.
$$
  The
singular frame $\te_0,\te_1, \te_2$ consists of the eigenvectors of the matrix \eqref{eqno:26},
namely
$$
\te_0=\frac1m \big(-\dt x_3\,e_2+\al x_2\,e_3\big),\quad \te_1=e_1,\quad \te_2=\frac1m \big(\al
x_2\,e_2+\dt x_3\,e_3\big)
$$
To find $\tf_1$ and $\tf_2$, calculate $A_\xi \te_1$ and $A_\xi\te_2$:
$$
A_\xi \te_1=\bt\big( x_3\,e_2-x_2\,e_3 \big),\quad A_\xi \te_2=-m \,e_1.
$$
Denote $\ve=sign(\,\bt)$. Then
$$
\tf_1=\ve \big( x_3\,e_2-x_2\,e_3 \big),\quad \tf_2=-e_1.
$$
Now we have
$$
\tilde
\Omega_{\sigma|00}=-\frac{1}{\sqrt{1+\lambda_\sigma^2}}\la(\nabla_{\te_0}A_\xi)\te_0,\tf_\sigma\ra.
$$
If $\xi$ is totally geodesic, then $\xi$ satisfy
$$
0=(\nabla_{\te_0}A_\xi)\te_0=\nabla_{\te_0}(A_\xi\te_0)-A_\xi\nabla_{\te_0}\te_0=A_\xi
A_{\te_0}\te_0
$$
Since \eqref{eqno:25} is applicable to any left-invariant unit vector field, we easily calculate
\begin{multline*}
A_{\te_0}\te_0=\frac{1}{m^2}
\left(%
\begin{array}{ccc}
 0  & \al\dt x_3 & -\dt\al x_2\ \\[1ex]
 \bt \al x_2 & 0 &  0 \\[1ex]
 \bt\dt x_3 & 0 & 0 \\
\end{array}%
\right)
\left(%
\begin{array}{c}
  0 \\[1ex]
  -\dt x_3 \\[1ex]
  \al x_2 \\
\end{array}%
\right) =\\[2ex]
 -\frac{1}{m^2}\al\dt(\dt x_3^2+\al x_2^2)\,\te_1.
\end{multline*}
Therefore,
$$
A_\xi A_{\te_0}\te_0=  -\frac{1}{m^2}\al\dt(\dt x_3^2+\al x_2^2) A_\xi \te_1=-\ve \bt\al\dt(\dt
x_3^2+\al x_2^2)\tf_1.
$$
Since $\bt\ne0,\al\ne0$ and $\dt\ne0$, we have
$$
\left\{
\begin{array}{c}
\al x_2^2+\dt x_3^2=0, \\[1ex]
 x_2^2+ x_3^2=1. \\
\end{array}
\right.
$$
Solving the system, we get
$$
x_2^2=\frac{-\dt}{\al-\dt},\quad x_3^2=\frac{\al}{\al-\dt}.
$$
Remind that $ \al+\dt>0,\quad \al\geq\dt$ by the choice of the frame. Therefore, the solution
exists, if $\dt<0$ and, as a consequence, $ \al>0$. Thus,
$$
\xi=  \pm\sqrt{\frac{-\dt}{\al-\dt}}\ e_2\pm\sqrt{\frac{\al}{\al-\dt}}\ e_3.
$$
Denote $\theta=\pm\,1$. Without loss of generality we can set
$$
\xi=  \theta\sqrt{\frac{-\dt}{\al-\dt}}\ e_2+\sqrt{\frac{\al}{\al-\dt}}\ e_3.
$$
As a consequence
$$
m=\sqrt{\vphantom{\big|}\al^2\frac{-\dt}{\al-\dt}+\dt^2\frac{\al}{\al-\dt}}=\sqrt{-\al\dt}.
$$

Moreover
$$
\begin{array}{l}
\ds\frac{\al}{m}x_2=\theta\,\frac{\al}{\sqrt{-\al\dt}}\sqrt{\frac{-\dt}{\al-\dt}}=
\theta\,x_3,\\[2ex]
\ds\frac{\dt}{m}x_3=\frac{\dt}{\sqrt{-\al\dt}}\sqrt{\frac{\al}{\al-\dt}}=
\frac{-\sqrt{(-\dt)^2}}{\sqrt{-\al\dt}}\sqrt{\frac{\al}{\al-\dt}}= -\theta\,x_2
\end{array}
$$
and we have
$$
\begin{array}{l}
\te_0=\frac{1}{m} \big(-\dt x_3\,e_2+\al
x_2\,e_3\big)=\theta\,\xi,\\[1ex]
\te_1=e_1=-\tf_2, \\[1ex]
\te_2=\frac{1}{m}\big(\al x_2\,e_2+\dt x_3\,e_3\big)= \theta(x_3\,e_2- x_2\,e_3)=\theta\ve\,\tf_1.
\end{array}
$$
With respect to this frame, we have
$$
\begin{array}{l}
A_\xi\te_0=A_\xi\xi=0,\\[1ex]
A_\xi\te_1=|\bt|\tf_1=\theta\ve|\bt|\,\te_2=\theta\bt\te_2\\[1ex]
A_\xi \te_2=m\tf_2=-m\,\te_1
\end{array}
$$
and the matrix $A_\xi$ takes the form
$$
A_\xi=
\left(%
\begin{array}{ccc}
  0 & 0 & 0 \\
  0 & 0 & -m \\
  0 & \theta\bt & 0 \\
\end{array}%
\right).
$$
A simple calculation  yields
\begin{equation}\label{eqno:27}
\begin{tabular}{|c|c|c|c|}
  \hline
  $\nabla$ & $\te_0$ & $\te_1$ & $\te_2$ \\
  \hline\hline
  \vphantom{\Big|}$\te_0$ & 0 & $- \theta m\te_2$ & $ \theta m\te_1$ \\
  \hline
  \vphantom{\Big|}$\te_1$ & $-\bt \te_2$ & 0 & $\bt\te_0$ \\
  \hline
  \vphantom{\Big|}$\te_2$ & $\theta m\te_1$ & $-\theta m\te_0-(\al+\dt)\te_2$ & $(\al+\dt)\te_1$ \\
  \hline
\end{tabular}\ .
\end{equation}
With respect to new frame, the derivatives $(\nabla_{\te_i}A_\xi)\te_k$ form the following Table.
\begin{center}
\begin{tabular}{|c|c|c|c|}
  \hline
  \vphantom{\Big|}     & $\te_0$ & $\te_1$ & $\te_2$ \\
  \hline\hline
  \vphantom{\Big|}$\te_0$ & 0 & $-m(\theta m-\bt)\,\te_1$ & $m(\theta m-\bt)\,\te_2$ \\
  \hline
  \vphantom{\Big|}$\te_1$ & $-m\bt\, \te_1$ & $\theta \bt^2 \,\te_0$ & 0 \\
  \hline
  \vphantom{\Big|}$\te_2$ & $-m\bt\,\te_2$ & $-(\al+\dt)( m-\theta\bt)\,\te_1$ & $\theta m^2\te_0+(\al+\dt)(m-\theta\bt)\,\te_2$ \\
  \hline
\end{tabular}\ .
\end{center}

Finally, the necessary components of the curvature operator can be found from the latter Table and
take the form
\begin{equation}\label{eqno:28}
\begin{array}{l}
\ds R(\te_0,\te_1)\xi=m(\theta m-2 \bt)\,\te_1,\\[1ex]
\ds R(\te_0,\te_2)\xi=-\theta m^2\,\te_2,\\[1ex]
\ds R(\te_1,\te_2)\xi=-(\al+\dt)(m- \theta\bt)\,\te_1.\\[1ex]
\end{array}
\end{equation}
Remark, also, that
$$
\tf_1=\theta\ve \te_2,\quad \tf_2=-\te_1.
$$
 Now, we can find all the entries of the
matrices $\tilde \Omega_\sigma$.
\begin{multline*}
\ds \tilde\Omega_{1|10}=\ds\frac{-\big< (\nabla_{\te_1}A_\xi)\te_0+  (\nabla_{\te_0}A_\xi)\te_1,
\tf_1
\big>+ \lambda_1^2 \big<R(\te_1, \te_0) \xi, \tf_1 \big> }{2(1+\lambda_1^2)}=\\
    \ds\frac{-\big< -\theta m^2\te_1, \theta\ve\te_2 \big>+
        \bt^2 \big<-m(\theta m-2\bt)\,\te_1, \theta\ve\te_2 \big>
        }{2(1+\bt^2)}  =0.
\end{multline*}

\begin{multline*}
\ds \tilde\Omega_{1|20}=\ds\frac{-\big< (\nabla_{\te_2}A_\xi)\te_0+ (\nabla_{\te_0}A_\xi)\te_2,
\tf_1 \big>+\lambda_1\lambda_2 \big<R(\te_1, \te_0) \xi, \tf_2
\big>}{2\sqrt{(1-\lambda_1^2)(1+\lambda_2^2)}} =\\[4ex]
 \ds \frac{-\big< m(\theta m-2\bt)\te_2, \theta\ve\te_2 \big>+
|\bt| m \big<-m(\theta m-2\bt)\te_1, -\te_1
\big>}{2\sqrt{(1+\bt^2)(1+m^2)}}= \\[3ex]
\ds \frac{m(m-2\theta \bt)(m|\bt|-\theta\ve)}{2\sqrt{(1-\bt^2)(1-m^2)}} = \frac{\ve m(\theta
m-2\bt)(m\bt-\theta)}{2\sqrt{(1+\bt^2)(1+m^2)}}=\\[1ex]
\frac{\ve\theta m(m-2\theta\bt)(m\bt-\theta)}{2\sqrt{(1+\bt^2)(1+m^2)}} .
\end{multline*}

\begin{multline*}
\ds \tilde\Omega_{1|11}=\frac{-\big< (\nabla_{\te_1}A_\xi)\te_1, \tf_1
\big>}{\sqrt{(1+\lambda_1^2)^3}} =\frac{-\big< \theta\bt^2\te_0,
\theta\ve\te_1 \big>}{(1+\bt^2)\sqrt{(1+\bt^2)}}=0.\\
\end{multline*}

\begin{multline*}
\ds \tilde\Omega_{1|12}=\ds\frac{-\big< (\nabla_{\te_1}A_\xi)\te_2+ (\nabla_{\te_2}A_\xi)\te_1,
\tf_1 \big>+\lambda_1^2\big<R(\te_1, \te_2) \xi,
\tf_1\big>}{2\sqrt{(1+\lambda_1^2)^2(1+\lambda_2^2)}} =\\[3ex]
\ds \frac{\big< (\al+\dt)(m-\theta  \bt)\,\te_1, \theta\ve \te_2\big>+\bt^2\big<-(\al+\dt)(m-\theta
\bt)\,\te_1,\theta\ve\te_2\big>}{2(1+\bt^2)\sqrt{1+m^2}} =0.
\end{multline*}

\begin{multline*}
\ds \tilde\Omega_{1|22}=\ds\frac{-\big< (\nabla_{\te_2}A_\xi)\te_2,\tf_1
\big>+\lambda_1\lambda_2\big<R(\te_1, \te_2) \xi,
\tf_2\big>}{\sqrt{(1+\lambda_1^2)(1+\lambda_2^2)^2}} =\\[3ex]
\! \frac{-\big<\theta   m^2\te_0+(\al+\dt)(m-\theta  \bt)\,\te_2, \theta\ve\te_2\big>+|\bt|
m\big<(\al+\dt)(m-\theta \bt)\,\te_1,\te_1\big>}{\sqrt{(1+\bt^2)}\,(1+m^2)}=\\[3ex]
\ds\frac{-\theta\ve(\al+\dt)( m-\theta \bt)+|\bt|m(\al+\dt)(m-\theta\bt)}{\sqrt{(1+\bt^2)}\,(1+m^2)}= \\[2ex]
 \ds\frac{\ve\theta(\al+\dt)(\theta m-\bt)(m\bt-\theta)}{\sqrt{(1+\bt^2)}\,(1+m^2)}.
\end{multline*}

Summarizing, we get
$$
\ds\tilde\Omega_1=
\left(%
\begin{array}{ccc}
  \ds0 &\ds 0 & \ds\frac12 \frac{\ve\theta m(m-2\theta\bt)( m\bt-\theta)}{\sqrt{(1+\bt^2)(1+m^2)}}
  \\[2ex]
  \ds0 &\ds 0 &\ds 0 \\[1ex]
\ds  \frac12 \frac{\ve\theta m(m-2\theta\bt)( m\bt-\theta)}{\sqrt{(1+\bt^2)(1+m^2)}} & \ds0&\ds \frac{\ve\theta(\al+\dt)(\theta m-\bt)( m\bt-\theta)}{\sqrt{(1+\bt^2)(1+m^2)}} \\
\end{array}%
\right).
$$
In a similar way, we find

\begin{multline*}
\ds\tilde\Omega_{2|10}=\ds\frac{-\big< (\nabla_{\te_1}A_\xi)\te_0+ (\nabla_{\te_0}A_\xi)\te_1,
\tf_2 \big>+ \lambda_2\lambda_1 \big<R(\te_2, \te_0) \xi, \tf_1 \big>}
{2\sqrt{(1+\lambda_2^2)(1+\lambda_1^2)}}=\\[3ex]
\ds\frac{-\big< -\theta m^2\te_1, -\te_1 \big>+ m\,|\bt| \big<\theta m^2\te_2, \theta\ve
\te_2\big>}{2\sqrt{(1+m^2)(1+\bt^2)}}=\ds\frac{ m^2 (m\bt-\theta)}{2\sqrt{(1+m^2)(1+\bt^2)}}.
\end{multline*}

\begin{multline*}
\ds\tilde\Omega_{2|20}=\ds\frac{-\big< (\nabla_{\te_2}A_\xi)\te_0+ (\nabla_{\te_0}A_\xi)\te_2,
\tf_2 \big>}{2\sqrt{(1+\lambda_2^2)^2}}=\ds\frac{-\big<
(m(\theta m-\bt)\te_2, -\te_1\big>} {2(1+m^2)}=0. \\
\end{multline*}

\begin{multline*}
\ds\tilde\Omega_{2|11}=\ds\frac{-\big< (\nabla_{\te_1}A_\xi)\te_1, \tf_2 \big>+ \lambda_2\lambda_1
\big<R(\te_2, \te_1) \xi, \tf_1 \big>}{\sqrt{(1+\lambda_2^2)(1+\lambda_1^2)}}=\\[3ex]
\ds\frac{-\big< \theta\bt^2\te_0, \te_1 \big>+ m |\bt| \big<(\al+\dt)(\theta m-\bt)\,\te_1,
\theta\ve\te_2\big>}{\sqrt{(1+m^2)(1+\bt^2)}}=0.
\end{multline*}

\begin{multline*}
\ds\tilde\Omega_{2|12}=\ds\frac{-\big< (\nabla_{\te_1}A_\xi)\te_2+ (\nabla_{\te_2}A_\xi)\te_1,
\tf_2 \big>+ \lambda_2^2\big<R(\te_2, \te_1) \xi, \tf_2 \big>
}{2\sqrt{(1+\lambda_2^2)^2(1+\lambda_1^2)}}=\\[3ex]
\ds\frac{\big< (\al+\dt)( m-\theta\bt)\,\te_1, -\te_1 \big>+
m^2\big<(\al+\dt)(m- \theta\bt)\,\te_1, -\te_1 \big>}{2(1+m^2)\sqrt{(1+\bt^2)}}=\\[3ex]
\ds \frac{-(\al+\dt)( m-\theta\bt)}{2\sqrt{1+\bt^2}}.
\end{multline*}

\begin{multline*}
\!\tilde\Omega_{2|22}=\ds\frac{-\big<
(\nabla_{\te_2}A_\xi)\te_2, \tf_2 \big> }{\sqrt{(1+\lambda_2^2)^3}}=
\ds\frac{-\big<\theta m^2\te_0+(\al+\dt)( m-\theta\bt)\,\te_2,
-\te_1\big> }{(1+m^2)\sqrt{1+m^2}}=0.\\
\end{multline*}
Summarizing, we get
$$
\!\!\tilde\Omega_2=
\!\left(%
\begin{array}{ccc}
  \ds0 & \ds \frac{m^2( m\bt-\theta)}{2\sqrt{(1+\bt^2)(1+m^2)}} & \ds0
  \\[2ex]
\ds \frac{m^2( m\bt-\theta)}{2\sqrt{(1+\bt^2)(1+m^2)}} &\ds 0 & \ds\frac{(\al+\dt)(\theta\bt- m)}{2\sqrt{(1+\bt^2)}} \\[1ex]
\ds0& \ds \frac{(\al+\dt)(\theta\bt- m)}{2\sqrt{(1+\bt^2)}}& 0 \\
\end{array}%
\!\right).
$$
Thus, for totally geodesic field $\xi$ we have a unique possible solution
$$
\bt=\theta m , \quad  m \bt=\theta.
$$
It follows, then,
$$
{-\al\dt}=m^2=1,\quad \bt=\theta.
$$
As a consequence,
$$
\pm\,\xi=\theta\frac{1}{\sqrt{1+\al^2}}\, e_2+\frac{\al}{\sqrt{1+\al^2}}
$$
is the corresponding totally geodesic unit vector field.
\end{proof}

\subsection{Geometrical description of totally geodesic unit vector field and the group}

\begin{proposition}\label{prop:3.1} Let $G$ be a non-unimodular three-dimensional Lie group with a left-invariant
metric. Suppose $G$ admits a left-invariant totally geodesic unit vector field $\xi$. Then either
\begin{itemize}
\item $G=L^2(-\alpha^2)\times E^1$, where $L^2(-\alpha^2)$ is the Lobachevski  plane of curvature
$-\al^2$, and $\xi$ is a parallel unit vector field on $G$ tangent to Euclidean factor, or
\item  $G$ admits the Sasakian structure; moreover, $G$ admits two hyperfoliations $\mathcal{L}_1, \mathcal{L}_2$ such that
\begin{itemize}
\item[(i)] the foliations $\mathcal{L}_1$ and $ \mathcal{L}_2$ are intrinsically flat, mutually orthogonal and has a constant extrinsic curvature,
\item[(ii)] one of them, say $\mathcal{L}_2$, is minimal,
\item[(iii)] the integral trajectories of the field $\xi$ are $\mathcal{L}_1\cap \mathcal{L}_2$.
\end{itemize}
\end{itemize}
\end{proposition}
\begin{proof} Suppose $\xi$ is as in the hypothesis. Consider the case
$\bt=\dt=0$ and $\xi=e_3$ of the Theorem \ref{Thrm:3.2}. The bracket operations take the form
$$
[e_1,e_2]=\al e_2, \quad [e_1,e_3]=0, \quad [e_2,e_3]=0
$$
and we conclude that the group admits  three integrable distributions, namely, $e_1\wedge e_2$,
$e_1\wedge e_3$ and $e_2\wedge e_3$. The Table of the Levi-Civita connection takes the form
\begin{center}
\begin{tabular}{|c|c|c|c|}
  \hline
  \vf $\nabla$ & $e_1$ & $e_2$ & $e_3$ \\[1ex]
  \hline
\vf  $e_1$ & 0 & $0$ & $ 0$ \\
  \hline
\vf  $ e_2 $ & $-\al e_2$ & $-\al e_1$ & $0$ \\
  \hline
\vf  $e_3$ & $ 0$ & $0$ & $0$ \\
  \hline
\end{tabular}
\end{center}
The only non-zero component of the curvature tensor of the group is of the form
$$
R(e_1,e_2)e_2=-\al^2 e_1.
$$
Thus, $G=L^2(-\al)\times R^1$ and the field $\xi=e_3$ is a parallel unit vector field on $G$
tangent to the Euclidean factor.

Consider the second case of the Theorem \ref{Thrm:3.2}.  If  $\bt=\theta ,m=\sqrt{-\al\dt}=1$, then
with respect to the singular frame the matrix $A_\xi$ takes the form
$$
A_\xi=
\left(%
\begin{array}{ccc}
  0 & 0 & 0 \\
  0 & 0 & -1 \\
  0 & 1 & 0 \\
\end{array}%
\right)
$$
and hence, $\xi=\theta\te_0$ is the Killing unit vector field. Therefore, by the Theorem
\ref{Thrm:2.3}, the structure
$$
\big(\phi=A_\xi, \ \xi,\ \eta=\la \xi,\cdot\ra\big);
$$
is Sasakian.

We can also say more about this Sasakian structure. The Table \eqref{eqno:27} in the case under
consideration takes the form
\begin{equation}\label{eqno:29}
\begin{tabular}{|c|c|c|c|}
  \hline
  $\nabla$ & $\te_0$ & $\te_1$ & $\te_2$ \\
  \hline\hline
  \vphantom{\Big|}$\te_0$ & 0 & $-\theta\te_2$ & $\theta\te_1$ \\
  \hline
  \vphantom{\Big|}$\te_1$ & $- \theta\te_2$ & 0 & $\theta\te_0$ \\
  \hline
  \vphantom{\Big|}$\te_2$ & $\theta \te_1$ & $-\theta \te_0-(\al+\dt)\te_2$ & $(\al+\dt)\te_1$ \\
  \hline
\end{tabular}
\end{equation}
an hence, for the  brackets we have
\begin{equation}\label{eqno:30}
[\te_0,\te_1]=0,\quad [\te_0,\te_2]=0,\quad [\te_1,\te_2]=2\theta \te_0+(\al +\dt)\,\te_2.
\end{equation}
From \eqref{eqno:30} we see that the distributions  $\te_0\wedge \te_2$ and $\te_0\wedge \te_1$ are
integrable. Denote by $\mathcal{L}_1$ and $\mathcal{L}_2$ the corresponding foliations generated by
these distributions. Then the integral trajectories of the field $\xi$ are exactly
$\mathcal{L}_1\cap \mathcal{L}_2$.

Denote $\Omega^{(1)}$ and $\Omega^{(2)}$ a second fundamental form of the $\mathcal{L}_1$ and
$\mathcal{L}_2$ respectively.  Since $\te_1$ and $\te_2$ are the unit normal vector fields for the
corresponding foliations, from \eqref{eqno:29} we can easily find
$$
\Omega^{(1)}=
\left(%
\begin{array}{cc}
  0 & 1 \\
  1 & \al+\dt \\
\end{array}%
\right), \qquad \Omega^{(2)}=
\left(%
\begin{array}{cc}
  0 & -\theta \\
  \theta & 0 \\
\end{array}%
\right)
$$
and see that $\mathcal{L}_2$ is a minimal foliation.

Setting $\xi=\theta \te_0$, we can find from \eqref{eqno:28} the corresponding curvature components
$$
R(\te_0,\te_2)\te_0=-\te_2, \quad  R(\te_0,\te_1)\te_0=-\te_1.
$$
Denote by $K^{(i)}_{int}$ and  $K^{(i)}_{ext}$ the intrinsic and extrinsic curvatures of the
corresponding foliations ($i=1,2$). Then $K^{(i)}_{ext}=\big<R(\te_0,\te_i)\te_i,\te_0\big>=1$. The
Gauss equation implies
$$
K^{(i)}_{int}=K^{(i)}_{ext}+\det \Omega^{(i)}=0.
$$
Therefore, both of the foliations are \emph{intrinsically flat} and have a \emph{constant extrinsic
curvature} $K^{(i)}_{ext}=1$.
\end{proof}
\section{Appendix}
The structure \eqref{eqno:20} appears also in a different setting. Remind that the unit tangent
bundle $T_1M^n$ is a hypersurface in $TM^n$ with a unit normal vector $\xi^v$ at each point
$(q,\xi)\in T_1M^n$. Define a unit vector field $\bar \xi$, a 1-form $\bar\eta$ and a $(1,1)$
tensor field $\bar\varphi$ on $T_1M^n$ by
$$
\bar\xi=-J\xi^v=\xi^h, \quad JX=\bar\varphi X+\bar\eta(X)\xi^v,
$$
where $J$ is a natural almost complex structure on $TM^n$, acting as
$$
JX^v=-X^h,\quad JX^h=X^v.
$$
The triple $(\bar\xi,\bar\eta,\bar\varphi)$ form a standard almost contact structure on $T_1M^n$
with Sasaki metric $g_S$. This structure is not almost contact \emph{metric} one. By taking
$$
\tilde\xi=2\bar\xi=2\xi^h, \quad \tilde\eta=\frac12\bar\eta,\quad \tilde\varphi=\bar\varphi,\quad
g_{cm}=\frac14g_S
$$
at each point $(q,\xi)\in T_1M^n$, we get the \emph{almost contact metric structure}
$(\tilde\xi,\tilde\eta,\tilde\varphi)$ on $(T_1M^n,g_{cm})$.

In a case of a general almost contact metric manifold $(\tilde M,\tilde
\xi,\tilde\eta,\tilde\varphi,\tilde g)$ a submanifold $N$ is called \emph{invariant} if
$\tilde\varphi(T_pN)\subset T_pN$ and \emph{anti-invariant} if $\tilde\varphi(T_pN)\subset
(T_pN)^\perp$ for every $p\in N$.

A unit vector field $\xi$ on a Riemannian manifold $(M^n,g)$ is called invariant (anti-invariant)
is the submanifold $\xi(M^n)\subset (T_1M^n,g_{cm})$ is invariant (anti-invariant). Recently, Binh
T.Q., Boeckx E. and Vanhecke L. have considered this kind of unit vector fields  and proved the
following Theorem \cite{Binh}.
\begin{theorem}
A unit vector field $\xi$ on $(M^n, g)$ is invariant if and only if
$$
\big( \phi=A_\xi, \, \xi\,,\,\eta=\big<\xi,\cdot\big>_g\,\big)
$$
is an almost contact structure on $M^n$. In particular, $\xi$ is a geodesic vector field on $M^n$
and $n=2m+1$.
\end{theorem}
Summarizing the results of this section, we come to the following conclusion.
\begin{proposition}
Every left invariant non-parallel totally geodesic unit vector field on a three-dimensional Lie
group $G$ with a left-invariant metric generates the invariant submanifold in $(T_1G, g_{cm})$.
\end{proposition}

\vspace{2ex}

\noindent
Alexander Yampolsky,\\
Department of Geometry,\\
Faculty of Mechanics and Mathematics,\\
Kharkiv National University,\\
Svobody Sq. 4,\\
61077, Kharkiv,\\
Ukraine.\\
e-mail: AlexYmp@gmail.com


\begin{thebibliography}{30}

\bibitem{Besse} Besse A. Manifolds all of whose geodesics are closed, Springer-Verlag, 1978.

\bibitem{Binh} Binh T.Q., Boeckx E. and Vanhecke L., \emph{Invariant and Anti-invariant Unit Vector Fields,} to
appear

\bibitem{BX-V1}
Boeckx~E., Vanhecke~L. \textit{Harmonic and minimal radial vector fields.} Acta Math. Hungar. 90
(2001), 317-331.

\bibitem{BX-V2}
Boeckx~E., Vanhecke~L. \textit{Harmonic and minimal vector fields on tangent and unit tangent
bundles.} Differential Geom. Appl. 13 (2000), 77-93.

\bibitem{BX-V4}
Boeckx~E., Vanhecke~L. \textit{Isoparametric functions and harmonic and minimal unit vector
fields}, Contemp. Math. 288 (2001), 20--31.

\bibitem{B-Ch-N} Brito F., Chacon P. and Naveira A.,\emph{ On the volume of vector fields on spaces of constant sectional
 curvature,} Comment. Math. Helv. 79 (2004), 300--316.

\bibitem{GM}Gil-Medrano~O. \textit{Relationship between volume and energy of unit vector fields}, Diff. Geom. Appl. 15
(2001), 137-152.


\bibitem{GM-GD-Vh} Gil-Medrano~O., Gonz\'alez-D\'avila J.C., Vanhecke L. \textit{Harmonic and minimal invariant unit
vector fields on homogeneous Riemannian manifolds}, Houston J. Math. 27 (2001), 377-409.


\bibitem{GM-LF}
Gil-Medrano~O.,Llinares-Fuster~E. \textit{Minimal unit vector fields}, T\^ohoku Math. J. 54 (2002),
71 -- 84.

\bibitem{GD-V1}
Gonz\'alez-D\'avila J.C., Vanhecke L. \textit{Examples of minimal unit vector fields.} Ann Global
Anal. Geom. 18 (2000), 385-404.

\bibitem{GD-V2}
Gonz\'alez-D\'avila J.C., Vanhecke L. \textit{Minimal and harmonic characteristic vector fields on
three-dimensional contact metric manifolds}. J. Geom. 72 (2001), 65-76.

\bibitem{G-Z}
Gluck~H., Ziller~W. \textit{On the volume of a unit vector field on the three-sphere}. Comm. Math.
Helv. 61 (1986), 177-192.

\bibitem{Mn} Milnor J. \emph{Curvatures of left invariant metrics on Lie groups,} Adv. in Math, 21
(1976), 293 -- 329.

\bibitem{TS-V1}
Tsukada K.,Vanhecke~L. \textit{Minimality and harmonicity for Hopf vector fields}, Illinois J.
Geom. 45 (2001), 441 -- 451.

\bibitem{TS-V2}
Tsukada K.,Vanhecke~L. \textit{Invariant minimal unit vector fields on Lie groups,} Period. Math.
Hungar. 40 (2000), 123-133.

\bibitem{TS-V3} Tsukada K.,Vanhecke~L. \textit{Minimal and harmonic vector fields on
$G_2(C^{m+2})$ and its dual space}. Monatsh. Math. 130 (2000), 143-154.

\bibitem{Weigm} Weigmink G., \textit{Total Bending of Vector Fields on Riemannian Manifolds,}
Math. Ann. \textbf{303} (1995) 325 -- 344.

\bibitem{Ym1} Yampolsky A. \textit{On the mean curvature  of a unit vector field,} Math. Publ.
Debrecen, 60/1-2 (2002), 131 -- 155.

\bibitem{Ym2} Yampolsky A. \textit{On the intrinsic geometry of a unit vector field,}
Comment. Math. Univ. Carol. 43/2 (2002), 299-317.

\bibitem{Ym3} Yampolsky A. \textit{A totally geodesic property of Hopf vector fields,} Acta Math. Hungar. 101/1-2
(2003), 93-112.

\bibitem{Ym4} Yampolsky A. \textit{On extrinsic geometry of unit normal vector field of Riemannian hyperfoliation,}
 Math. Publ. Debrecen {63}/4 (2003), 555 -- 567.

\bibitem{Ym5}  Yampolsky A. {\em Full description of totally geodesic unit vector field on Riemannian
2-manifold}, Matematicheskaya fizika, analiz, geometriya 2004, 11/3, 355-365.

\bibitem{Ym6} Yampolsky A. \emph{On special types of minimal and totally geodesic unit vector
fields}, Proceedings of the  Sevent International Conference in Geometry, Integrability and
Quantization June 2 -- 10, 2005, Varna, Bulgaria (to appear).

\end{thebibliography}
\end{document}